\title{The Local Distribution of the Number of Small Prime Factors - Variation of the Classical Theme}
\author{Krishnaswami Alladi and Todd Molnar}
\documentclass{article}
\usepackage{graphicx}
\usepackage{amsmath}
\usepackage{amsfonts}
\usepackage{amssymb}

\newenvironment{proof}{{\sc Proof:}}{~\hfill QED}

\begin{document}
\newpage
\maketitle

ABSTRACT: We obtain uniform estimates for $N_k(x,y)$, the number of
positive integers $n$ up to $x$ for which $\omega_y(n)=k$, where $\omega_y(n)$
is the number of distinct prime factors of $n$ which are $<y$. The motivation
for this problem is an observation due to the first author in 2015 that for
certain ranges of $y$, the asymptotic behavior of $N_k(x,y)$ is different from
the classical situation concerning $N_k(x,x)$ studied by Sathe and Selberg.
We demonstrate this variation of the classical theme; 
to estimate $N_k(x,y)$ we study the sum $S_z(x,y)=\sum_{n\le x}z^{\omega_y(n)}$ for
$Re(z)>0$ by the Buchstab-de Bruijn method. We also
utilize a certain recent result of Tenenbaum to complete our asymptotic
analysis.    


{\textit{2000 Subject Classification:}} 11M06, 11M41, 11N25, 11N37, 11N60

\section{Introduction}

\indent The function
$$
\omega_{y}(n)=\sum_{\substack{p|n\\
		p<y}}1,
$$
\noindent $p$ a prime number, figures prominently in the proofs of the Erd{\H o}s-Kac theorem which concerns the global distribution of the number of prime factors.  Here we focus on the ``local distribution" of $\omega_{y}(n)$, that is, we study the function
$$
N_{k}(x,y)=\sum_{\substack{n\leq x\\
		\omega_{y}(n)=k}}1
$$
\noindent with emphasis on results which are uniform.  When $\alpha=\log(x)/\log(y)>1$, and $k$ is small, the behavior of $N_{k}(x,y)$ is different from the classical case of $N_{k}(x)=N_{k}(x,x)$ (see Theorems 2, 3, 11, and 12), a phenomenon first observed by Alladi in 2015.  We will supply the analysis necessary to explain this phenomenon.  However, as $k$ approaches $\log\log(y)$, or as $y$ approaches $x$, the behavior is similar to the classical situation.  We investigate how and when such a transition takes place.

\indent In preparation for this investigation we shall study the behavior of the sum
\begin{equation}
S_{z}(x,y)=\sum_{n\leq x}z^{\omega_{y}(n)}
\end{equation}
\noindent where $z\in\mathbb{C}$. Note that in the special case $z=0$,
$S_0(x,y)=\Phi(x,y)$, the well known function counting the number of
uncancelled elements in the Sieve of Eratosthenes. Obtaining asymptotic
estimates for $S_z(x,y)$ which are both sharp enough for application and
uniform in $y$ will require a variety of tools. Besides the elementary
techniques of Section 3, we shall also employ analytic methods in Sections 2,
4, and 5 such as de Bruijn's method of utilizing Buchstab iterations and
difference-differential equations [5].  It will be important to note that
the ranges for the estimates in these sections overlap, thereby permitting
a result valid in a larger range.  More importantly, by comparing the
estimates in these two ranges, we obtain an asymptotic estimate for a certain
function satisfying a difference-differential equation without going through
a saddle point analysis; this technique was used by Alladi in [4].

\indent Once we have asymptotic estimates for $S_{z}(x,y)$, we may then study $N_{k}(x,y)$ by recognizing that
\begin{equation}
N_{k}(x,y)=\frac{1}{2\pi i}\int_{\gamma}\frac{S_{z}(x,y)}{z^{k+1}}dz
\end{equation}
\noindent for a suitable contour $\gamma$.  The idea of using such a contour
integral is due to Selberg [10] who studied the case $y=x$.  If we let
$$
g(s,z)=\prod_{p}\left(1+\frac{z}{p^{s}-1}\right)\left(1-\frac{1}{p^{s}}\right)^{z}
$$ 
\noindent then Selberg demonstrated that
$$
S_{z}(x)=\frac{g(1,z)}{\Gamma(z)}\frac{x}{\log^{1-z}(x)}+O\left(\frac{x}{\log^{2-z}(x)}\right)
$$
\noindent uniformly for $|z|\leq R$, where $\Gamma(z)$ is the Gamma function.  Then using the integral in (2) with $\gamma$ a circle centered at the origin with radius 
$$
\rho=\frac{k-1}{\log\log(x)},
$$ 

\noindent Selberg demonstrated that
\begin{equation}
N_{k}(x)=\frac{x}{\log(x)}\frac{g\left(1,\frac{k-1}{\log\log(x)}\right)}{\Gamma\left(1+\frac{k-1}{\log\log(x)}\right)}\frac{(\log\log(x))^{k-1}}{(k-1)!}\left(1+O\left(\frac{k}{(\log\log(x))^{2}}\right)\right)
\end{equation}
\noindent uniformly for $k\leq R\log\log(x)$, where $R$ is any fixed positive
number.  This result improves upon that of Landau (see(7) below) and
Sathe in [9].  Our methods leading to Theorems 10 and 12 show that for certain
ranges of $k$ and $y$ 
$$
N_{k}(x,y)\asymp N_{k+1}(x);
$$
\noindent in such a situation, the function $N_{k}(x,y)$ is of the size of
the $(k+1)$-st
Landau function.  The above equation is the precise
statement of the phenomenon observed by Alladi in 2015, and our results
pertaining to this were established in 2016. 

\indent In a fundamental paper Hal\'{a}sz [7] studied the local distribution
of the general additive function
$$
\omega(n;E):=\sum_{\substack{p|n\\
		p\in E}}1
$$

\noindent where $E$ is an arbitrary set of prime numbers.  He obtained
asymptotic estimates for
$$
N_{k}(x;E):=\sum_{\substack{n\leq x\\
		\omega(n;E)=k}}1
$$
\noindent when $\omega(n;E)\sim E(x)$, where
$$
E(x)=\sum_{\substack{p\leq x\\
		p\in E}}\frac{1}{p}
$$
\noindent (see Elliott's book [6] for the proof). Note here that
$\omega(n;E)$ will be almost always the size of $E(x)$. In a recent paper
Tenenbaum [13] has strengthened Hal\'{a}sz's result by
extending the range of the asymptotic formula to 
$$
1/\kappa\leq k/E(x)\leq\kappa
$$

\noindent for any large $\kappa>0$; his results are uniform in that range.
In the course of explaining the phenomenon observed by Alladi, Tenenbaum [14]
communicated to us in 2016 that by choosing $E=\{p|p<y\}$,
one may use the Selberg-Delange method and other techniques to obtain
sharper error terms than ours in some instances.  The emphasis of this paper
is to obtain sufficiently
sharp and uniform asymptotic estimates that will demonstrate and explain
the phenomenon observed by Alladi, which seems to have escaped attention.
Our methods are different from those of Hal\'{a}sz and Tenenbaum and are of
intrinsic interest. 

\section{Asymptotic estimate of $S_{z}(x,y)$ for small $y$}

\indent Following the method of de Bruijn in [5] we may obtain an estimate for $S_{z}(x,y)$ using contour integration.  The representation is valid for any complex $z\in\mathbb{C}$ and will yield an asymptotic estimate provided $z\neq-p+1$, with $p$ a prime number.  First note that as $z^{\omega_{y}(n)}$ is multiplicative we obtain
$$
F(s):=\sum_{n=1}^{\infty}\frac{z^{\omega_{y}(n)}}{n^{s}}=\zeta(s)\prod_{p<y}\left(1+\frac{z-1}{p^{s}}\right)=\zeta(s)g(s,y,z)
$$
\noindent where $\zeta(s)$ is the Riemann zeta-function and
$$
g(s,y,z)=\prod_{p<y}\left(1+\frac{z-1}{p^{s}}\right),
$$
\noindent so that the function $g(s,y,z)$ is now analytic in both $s$ and $z$.

\bigskip

\noindent\textbf{Lemma 1:} \textit{If $\Re(s)=\sigma\geq1$ and $|z|\leq R$ then
$$
g(s,y,z)<<_{R}\log^{|z-1|}(y);
$$
\noindent if $b:=\max\left(1-\frac{1}{\log(T)},1-\frac{1}{\log(y)}\right)$ and $\sigma\geq b$ then there exists a constant $C\in\mathbb{R}^{+}$ such that
$$
g(s,y,z)<<_{R}\log^{C|z-1|}(y).
$$
}

\bigskip

\begin{proof} The proof follows using the methods in [5] and standard estimates, and is omitted.

\end{proof}

\noindent With the estimate of Lemma 1 we may now prove

\bigskip

\noindent\textbf{Theorem 1:} \textit{Let $R>0$ be fixed, $|z|\leq R$, and $x\geq y\geq3$ then
	$$
	S_{z}(x,y)=x\prod_{p<y}\left(1+\frac{z-1}{p}\right)+O(xe^{-\alpha}\log^{D}(x))+O_{R}\left(\frac{x}{\log^{R+2}(x)}\right)
	$$
	\noindent for some absolute positive constant $D>0$.}
	
\bigskip

\begin{proof} From a standard inversion procedure (the effective Perron integral formula in [11]) we have the following integral representation
$$
S_{z}(x,y)=\sum_{n\leq x}z^{\omega_{y}(n)}
$$
\begin{equation}
=\frac{1}{2\pi i}\int_{a-iT}^{a+iT}\zeta(s)g(s,y,z)\frac{x^{s}}{s}ds+O\left(\sum_{n=1}^{\infty}\left(\frac{x}{n}\right)^{a}\tau_{R}(n)\min\left(1,\frac{1}{T|\log(x/n)|}\right)\right).
\end{equation}

\noindent The sum in the above equation (4) is estimated using a standard technique given explicitly in [11] (although some variation of this technique is used in almost every application of the Perron integral formula); dividing the sum into three intervals: $n<x/2$, $n>3x/2$, and $x/2\leq n\leq3x/2$ for $n\neq x$ we find that
$$
\sum_{n=1}^{\infty}\left(\frac{x}{n}\right)^{a}\tau_{R}(n)\min\left(1,\frac{1}{T|\log(x/n)|}\right)<<_{R}\frac{x}{\log^{R+2}(x)}+\frac{x\log^{2R}(x)}{T}.
$$

\noindent We omit the details which are standard.  Choosing $T:=e^{\sqrt{\log(x)}}$ in the above we obtain

\begin{equation}
S_{z}(x,y)-\frac{1}{2\pi i}\int_{a-iT}^{a+iT}\zeta(s)g(s,y,z)\frac{x^{s}}{s}ds<<_{R}\frac{x}{\log^{R+2}(x)}.
\end{equation}

\noindent It therefore remains to evaluate the integral in equation (5).

\indent Let $\Gamma$ be the rectangular contour with vertices $a+iT, b+iT, b-iT,$ and $a-iT$ with
$$
\max\left(1-\frac{1}{\log(T)},1-\frac{1}{\log(y)}\right)=b<1<a=1+\frac{1}{\log(x)}
$$
\noindent and $T=e^{\sqrt{\log(x)}}$ defined as above.  Using Cauchy's theorem and taking into account the simple pole of $\zeta(s)$ with residue $1$ along with some standard estimates of the vertical and horizontal segments of the contour we obtain (with the aid of Lemma 1) the statement of the theorem.  The two error terms arise from the two estimates in the definition of $b$. 

\end{proof}

\indent Theorem 1 can of course be used to obtain an effective asymptotic estimate for $S_{z}(x,y)$ for certain ranges of $y$.  In fact, if $\alpha$ is large enough to suppress the second term in Theorem 1 then, provided $y\rightarrow\infty$ and $z\neq1$ or $1-p$, the main term will simply be
$$
x\prod_{p<y}\left(1+\frac{z-1}{p}\right)\asymp_{z}\frac{x}{\log^{1-z}(y)}.
$$
\noindent This can be improved with the following

\bigskip

\noindent\textbf{Corollary 1:} \textit{Let $z\in\mathbb{C}$, $|z|\leq R$, and let $D$ denote the constant in Theorem 1.  If $\alpha\geq(R+D+1+\epsilon)\log\log(x)$, $\epsilon>0$ an arbitrary small fixed constant, then
	$$
	S_{z}(x,y)=x\prod_{p<y}\left(1+\frac{z-1}{p}\right)+O\left(\frac{x}{\log^{R+1+\epsilon}(y)}\right);
	$$
	\noindent in particular, there exists a positive constant $K$ such that if $\alpha\geq K\log\log(x)$ and $z\neq1,$ or $1-p$ then}
	$$
	S_{z}(x,y)\sim x\prod_{p<y}\left(1+\frac{z-1}{p}\right).
	$$

\bigskip

\indent Following de Bruijn in [5] we may introduce an additional parameter
$\lambda>0$ by replacing $1-\frac{1}{log(y)}$ with $1-\frac{\lambda}{\log(y)}$,
in defining $b$. Then carrying out the analysis of Theorem 1, we may optimize
this parameter to achieve an error term which is
$<<xe^{\alpha\log\alpha}\log^{D}(x)$.  We do not require the full strength of this
result and so we settled for the simpler bound in Theorem 1.  Theorem 1 and Corollary 1 can also be obtained using the more elementary technique of convolution sums, as in Hall and Tenenbaum [8].  Note that if $\alpha$ satisfies the
range of Corollary 1, then
\begin{equation}
y=x^{1/\alpha}\leq e^{\frac{\log(x)}{K\log\log(x)}}.
\end{equation}
\noindent This range of $y$ will be crucial for subsequent results. 

\section{Estimate of $N_{k}(x,y)$ for very large $y$}

\indent In this section we will derive estimates for $S_{z}(x,y)$ when $y$ is \textit{very large}, by which we mean that if $\beta=x/y$, then $\beta$ is fixed.  It is a classical result due to Edmund Landau that if $k$ is fixed and $\beta=1$ then
\begin{equation}
N_{k}(x)=\frac{x(\log\log(x))^{k-1}}{(k-1)!\log(x)}+O\left(\frac{x(\log\log(x))^{k-2}}{(k-2)!\log(x)}\right).
\end{equation}
\noindent This asymptotic is implied by Selberg's result of (3). We shall see in Theorem 2 that the asymptotic behavior of $N_{k}(x,y)$ is quite different from (6) for certain ranges of $y$.

\bigskip

\noindent\textbf{Lemma 2:} \textit{For any fixed integer $k>1$ and $\ell\geq1$ we have}
$$
\displaystyle\sum_{p_{1}^{e_{1}}p_{2}^{e_{2}}...p_{k}^{e_{k}}<x}\frac{\log^{\ell}(p_{1}^{e_{1}}...p_{k}^{e_{k}})}{p_{1}^{e_{1}}p_{2}^{e_{2}}...p_{k}^{e_{k}}}=\frac{\log^{\ell}(x)(\log\log(x))^{k-1}}{(k-1)!\ell}+O\left(\frac{\log^{\ell}(x)(\log\log(x))^{k-2}}{(k-2)!\ell}\right).
$$

\bigskip

\begin{proof}
Let $k,\ell\in\mathbb{Z}$ be fixed with $k>1$ and $\ell\geq1$.  By Stieltjes integration
$$
\displaystyle\sum_{p_{1}^{e_{1}}p_{2}^{e_{2}}...p_{k}^{e_{k}}<x}\frac{\log^{\ell}(p_{1}^{e_{1}}...p_{k}^{e_{k}})}{p_{1}^{e_{1}}p_{2}^{e_{2}}...p_{k}^{e_{k}}}=\int_{2}^{x}\frac{\log^{\ell}(t)}{t}dN_{k}(t),
$$
\noindent and applying integration by parts, we get
$$
\displaystyle\int_{2}^{x}\frac{\log^{\ell}(t)}{t}dN_{k}(t)=N_{k}(x)\frac{\log^{\ell}(x)}{x}-\int_{2}^{x}N_{k}(t)\left(\frac{\ell\log^{\ell-1}(t)-\log^{\ell}(t)}{t^{2}}\right)dt.
$$

\noindent From the above and (7) we may obtain the statement of the lemma
by one more intergration.

\end{proof}

\newpage

\noindent\textbf{Remarks:}

\noindent (i) Even though we have established Lemma 2 for $\ell\geq1$ and $k\geq2$ we will use it in the sequel only for $\ell=1$.

\noindent (ii) When $k=\ell=1$ there is the well-known result
$$
\sum_{p^{e}\leq x}\frac{\log(p)}{p^{e}}=\log(x)+c+O\left(\frac{1}{\log(x)}\right).
$$

\noindent (iii) The estimate given in Lemma 3 still applies if the summation is over squarefree integers, that is, if $\ell\geq1$ and $k>1$
\begin{equation}
\sum_{p_{1}...p_{k}<y}\frac{\log^{\ell}(p_{1}...p_{k})}{p_{1}...p_{k}}=\frac{\log^{\ell}(x)(\log\log(x))^{k-1}}{(k-1)!\ell}+O\left(\frac{\log^{\ell}(x)(\log\log(x))^{k-2}}{(k-2)!\ell}\right).
\end{equation}
\noindent To establish (8) we follow the proof given above for Lemma 2 except instead of the integrator $N_{k}(x)$ we use 
$$
N_{k}^{\ast}(x)=\sum_{\substack{n\leq x\\
		\Omega(n)=\omega(n)=k}}1
$$
\noindent which is the number of square-free integers with precisely $k$ prime factors.  The estimates of Lemma 2 and equation (8) agree because for fixed $k$
$$
N_{k}(x)= N_{k}^{\ast}(x)+O\left(\frac{x(\log\log(x))^{k-2}}{\log(x)(k-2)!}\right).
$$

\indent Let $\beta:=x/y$ be fixed.  For all but a finite number of integers $x>\beta^{2}$.  We begin by considering $N_{k}(x,y)$ when $k=0$, that is
\begin{equation}
N_{0}(x,y)=\sum_{\substack{n\leq x\\
		\omega_{y}(n)=0}}1=\sum_{\substack{n\leq x\\
		p^{-}(n)\geq y>\sqrt{x}}}1=\Phi(x,y).
\end{equation}

\noindent If $y>\sqrt{x}$ then the sum in (9) will be counting integers $n\leq x$ with all prime factors $>\sqrt{x}$; however, the only integers with this property are $1$ and the prime numbers $y\leq p\leq x$.  Hence, we have 

\bigskip

\noindent\textbf{Lemma 3:} \textit{If $x\geq y>\sqrt{x}$ then}
$$
N_{0}(x,y)=1+\pi(x)-\pi(y-1).
$$

\bigskip

\noindent From the prime number theorem we may immediately obtain an asymptotic estimate in the form
$$
N_{0}(x,y)=\frac{x}{\log(x)}+O\left(\frac{x}{\log^{2}(x)}\right)-\frac{y}{\log(y)}+O\left(\frac{y}{\log^{2}(y)}\right)
$$
\begin{equation}
=\frac{x}{\log(x)}-\frac{y}{\log(y)}+O\left(\frac{x}{\log^{2}(x)}\right),
\end{equation}
\noindent with obvious improvements being possible.  Comparing this with $N_{0}(x,x)=1$ (the corresponding Landau term) we can immediately see that $N_{0}(x,y)$ is much larger when $\beta>1$.

\indent Next, consider $N_{k}(x,y)$ when $k=1$.  The numbers which contribute to this sum will be of the form $n=mp^{e}\leq x$ with $p<y$, $p^{-}(m)\geq y>\sqrt{x}$, and $m\leq x/p^{e}$, hence,
$$
N_{1}(x,y)=\sum_{\substack{n\leq x\\
		\omega_{y}(n)=1}}1=\sum_{\substack{mp^{e}\leq x\\
		p<y\\
		p^{-}(m)\geq y>\sqrt{x}}}1=\sum_{\substack{p^{e}\leq x\\
		p<y}}1+\sum_{\substack{m\leq x/p^{e}\\
		p<y\\
		p^{-}(m)\geq y}}1
$$

\noindent where the first term corresponds to $m=1$ and the second to $m>1$; however, since $p^{-}(m)\geq y>\sqrt{x}$ this forces $m=q$ to be prime.  Consequently
$$
N_{1}(x,y)=\pi(y-1)+O(\sqrt{y})+\sum_{\substack{y\leq q\leq x/p^{e}\\
		p<y}}1.
$$
\noindent By evaluating the sum on the right in the above expression, we get

\bigskip

\noindent\textbf{Lemma 4:} \textit{If $x \geq y>1$ with $\beta=x/y$ fixed, then}
$$
N_{1}(x,y)=\frac{x}{\log(y)}\sum_{p^{e}<\beta}\left(\frac{1}{p^{e}}-\frac{1}{\beta}\right)+\frac{x}{\beta\log(y)}+O_{\beta}\left(\frac{x}{\log^{2}(x)}\right).
$$

\noindent\textbf{Remark:} In Lemma 4, the expression $\frac{x}{\beta\log(y)}$ is the Landau term (the result in (6)).  Note that it decreases as $\beta$ increases.

\indent If we now allow $\beta\rightarrow\infty$ slowly the estimate in lemma 4 implies that $N_{1}(x,y)\sim\frac{x\log\log\beta}{\log(x)}$, as this term dominates the Landau term of $\frac{x}{\beta\log(y)}$; meanwhile, if we allow $\beta\rightarrow1$ (equivalently $y\rightarrow x$) then $\sum_{p^{e}<\beta}\left(\frac{1}{p}-\frac{1}{\beta}\right)\rightarrow0$ (in fact, the sum is zero as soon as $\beta<2$) giving the classical result of
$$
N_{1}(x,x)=\frac{x}{\log(x)}+O\left(\frac{x}{\log^{2}(x)}\right).
$$

\indent We can apply the above techniques (and Lemma 2) to analyze the sum $N_{k}(x,y)$ for general fixed $k$ and $\beta>1$.

\bigskip

\noindent\textbf{Theorem 2:} \textit{For arbitrary fixed $\beta:=x/y>1$ and $k\in\mathbb{Z}$, $k>1$ fixed}
$$
N_{k}(x,y)=\frac{x}{\log(y)}\sum_{p_{1}^{e_{1}}...p_{k}^{e_{k}}\leq\beta}\left(\frac{1}{p_{1}^{e_{1}}...p_{k}^{e_{k}}}-\frac{1}{\beta}\right)+\frac{x(\log\log(y))^{k-1}}{\log(y)(k-1)!}
$$
$$
+O\left(\frac{x}{\log(y)}\frac{(\log\log(x))^{k-2}}{(k-2)!}\right).
$$

\bigskip

\noindent Theorem 2 is proved by noting that

$$
N_{k}(x,y)=\sum_{\substack{p_{1}^{e_{1}}...p_{k}^{e_{k}}\leq \beta\\
			e_{j}>0\\
			p_{j}<y}}\Phi\left(\frac{x}{p_{1}^{e_{1}}...p_{k}^{e_{k}}},y\right)+\sum_{\substack{\beta<p_{1}^{e_{1}}...p_{k}^{e_{k}}\leq x\\
			e_{j}>0\\
			p_{j}<y}}1,
$$

\noindent and utilizing the above lemmas.

\indent As a comment, if we now permit $\beta\rightarrow\infty$ slowly
then the first term in Theorem 2 will begin to dominate, which is
the phenomenon noticed in 2015.

\indent The above analysis will permit us to estimate $N_{k}(x,y)$ for all $y\in(\sqrt{x},x]$, provided $k$ is fixed, because the only fact used in the course of the proof of Theorem 2 is that if $x>x/n\geq y>\sqrt{x}$ then $\Phi(x/n,y)=\pi(x/n)-\pi(y-1)$.  In addition, the sum
$$
\sum_{p\leq\beta}\pi\left(\frac{x}{p}\right)\sim\sum_{p\leq\beta}\frac{x}{p\log(x/p)}=\frac{x}{\log(x)}\sum_{p\leq\beta}\frac{1}{p\left(1-\frac{\log(p)}{\log(x)}\right)}
$$
$$
=\frac{x}{\log(x)}\left(\sum_{p\leq\beta}\frac{1}{p}\right)+O\left(\frac{x}{\log^{2}(x)}\sum_{p\leq\beta}\frac{\log(p)}{p}\right)
$$
$$
=\frac{x}{\log(x)}\log\log\beta+O\left(\frac{x}{\log(x)}\right).
$$
\noindent This analysis can be continued when we have $k$ prime factors as
$$
\sum_{p_{1}...p_{k}\leq\beta}\pi\left(\frac{x}{p_{1}...p_{k}}\right)=\frac{x}{\log(x)}\sum_{p_{1}...p_{k}\leq\beta}\frac{1}{p_{1}...p_{k}}
$$
$$
+O\left(\frac{x}{\log^{2}(x)}\sum_{p_{1}...p_{k}}\frac{\log(p_{1}...p_{k})}{p_{1}...p_{k}}+\sum_{p_{1}...p_{k}\leq\beta}\frac{x}{p_{1}...p_{k}\log^{2}(x/p_{1}...p_{k})}\right),
$$
\noindent So that, with the aid of Lemma 2, we may arrive at the following

\bigskip

\noindent\textbf{Theorem 3:} \textit{For $\sqrt{x}<y\leq x$ we have}
$$
N_{k}(x,y)=\frac{x}{\log(x)}\frac{(\log\log\beta^{\ast})^{k}}{k!}+\frac{x}{\log(x)}\frac{(\log\log(x))^{k-1}}{(k-1)!}
$$
$$
+O\left(\frac{x(\log\log(x))^{k-2}}{\log(x)(k-2)!}+\frac{x}{\log(x)}\frac{(\log\log\beta^{\ast})^{k-1}}{(k-1)!}\right),
$$
\noindent \textit{where} $\beta^{\ast}=\max(\beta,10)$.

\bigskip

\indent We may rewrite Theorem 3 as

\begin{equation}
N_{k}(x,y)=w(\alpha)\frac{x}{\log(y)}\frac{(\log\log\beta^{\ast})^{k}}{k!}+\frac{x}{\log(x)}\frac{(\log\log(x))^{k-1}}{(k-1)!}
$$
$$
+O\left(\frac{x(\log\log(x))^{k-2}}{\log(x)(k-2)!}+\frac{x}{\log(x)}\frac{(\log\log\beta^{\ast})^{k-1}}{(k-1)!}\right),
\end{equation}

\noindent where $w(\alpha)=\frac{1}{\alpha}$ for $1\leq\alpha<2$.

We close this section with the following interesting consequence of the preceding theorem.

\bigskip

\noindent\textbf{Corollary 2:} \textit{Let $k\in\mathbb{Z}^{+}$ be fixed and $\sqrt{x}< y<xe^{-e^{(1+\epsilon)k^{1/k}(\log\log(x))^{1-1/k}}}$ for all $\epsilon>0$.  Then}
$$
N_{k}(x,y)\sim w(\alpha)\frac{x(\log\log(y))^{k}}{k!\log(y)}.
$$

\bigskip

\begin{proof} The corollary follows by comparing the first two terms of
equation (11), and noting that for the stated range the first term
dominates.
	
\end{proof}

\indent It is also possible to show that the results of Corollary 2 hold for
arbitrary fixed $\alpha>2$ and $k$ fixed.  For $y<\sqrt{x}$, we have for
$k=0$
$$
N_{0}(x,y)=\Phi(x,y)=w(\alpha)\frac{x}{\log(y)}+O\left(\frac{x}{\log^{2}(y)}\right),
$$
\noindent where
$$
w(\alpha)=\begin{cases}
\frac{1}{\alpha}, \text{ for } 1\leq\alpha\leq2\\
\frac{d}{d\alpha}\left(\alpha w(\alpha)\right)=w(\alpha-1), \text{ for } \alpha\geq2.
\end{cases}
$$ 
\noindent This is a well-known result for the number of uncanceled
elements in the sieve of Eratosthenes, and $w(\alpha)$ is known as
the Buchstab function.

\indent Let $\wp_{y}=\prod_{p<y}p$ and $k=1$, we see
$$
N_{1}(x,y)=\sum_{\substack{n=mp^{e}\leq x\\
		(m,\wp_{y})=1\\
		p<y}}1=\sum_{\substack{p<y\\
		e\geq1}}\Phi\left(\frac{x}{p^{e}},y\right)
$$
\begin{equation}
=\sum_{p<y}\Phi\left(\frac{x}{p},y\right)+\sum_{\substack{p<y\\
		e\geq2}}\Phi\left(\frac{x}{p^{e}},y\right)=\sum_{1}+\sum_{2}.
\end{equation}
\noindent Clearly,
\begin{equation}
\sum_{2}=\sum_{\substack{p^{e}\leq x/y\\
		p<y\\
		e\geq2}}\Phi\left(\frac{x}{p^{e}},y\right)+\sum_{\substack{p^{e}>x/y\\
		p<y,p^{e}\leq x\\
		e\geq2}}\Phi\left(\frac{x}{p^{e}},y\right).
\end{equation}
\noindent Using the estimate
$$
\Phi\left(x,y\right)<<\max\left(\frac{x}{\log(y)},1\right),
$$
\noindent we see that (13) is
\begin{equation}
<<\sum_{\substack{p<y\\
		p^{e}\leq x\\
		e\geq2}}\max\left(\frac{x}{p^{e}\log(y)},1\right)<<\frac{x}{\log(y)}+\sqrt{x}<<\frac{x}{\log(y)}.
\end{equation}
\noindent Regarding the first sum in (12), we use the fact that for $0\leq u\leq1$
$$
w(\alpha-u)=w(\alpha)+O(u),
$$
\noindent so that
$$
\sum_{1}=\sum_{p<y}\frac{x}{p}\frac{w\left(\frac{\log(x)-\log(p)}{\log(y)}\right)}{\log(y)}+O\left(\frac{x}{\log^{2}(y)}\sum_{p<y}\frac{1}{p}\right)
$$
$$
=w(\alpha)\frac{x}{\log(y)}\sum_{p<y}\frac{1}{p}+O\left(\frac{x}{\log^{2}(y)}\sum_{p<y}\frac{\log(p)}{p}\right)
$$
\begin{equation}
=w(\alpha)\frac{x\log\log(y)}{\log(y)}+O\left(\frac{x}{\log(y)}\right).
\end{equation}
\noindent Using (14) and (15) in (12) we conclude
$$
N_{1}(x,y)=w(\alpha)\frac{x\log\log(y)}{\log(y)}+O\left(\frac{x}{\log(y)}\right).
$$
\indent Similarly, we may also estimate $N_{2}(x,y)$ for fixed $\alpha>2$,
but this is a bit more complicated.  Note that
$$
N_{2}(x,y)=\sum_{\substack{n=mp_{1}^{e_{1}}p_{2}^{e_{2}}\leq x\\
		(m,\wp_{y})=1\\
		p_{1},p_{2}<y}}1=\sum_{p_{1},p_{2}<y}\Phi\left(\frac{x}{p_{1}^{e_{1}}p_{2}^{e_{2}}},y\right)
$$
\begin{equation}
=\sum_{p_{1},p_{2}<y}\Phi\left(\frac{x}{p_{1}p_{2}},y\right)+\sum_{\substack{p_{1},p_{2}<y\\
		e_{1}+e_{2}\geq3}}\Phi\left(\frac{x}{p_{1}^{e_{1}}p_{2}^{e_{2}}}\right)=\sum_{3}+\sum_{4}.
\end{equation}
\noindent In this case
$$
\sum_{4}<<\sum_{\substack{p_{1},p_{2}<y\\
		x/y\leq p_{1}^{e_{1}}p_{2}^{e_{2}}\leq x\\
		e_{1}+e_{2}\geq3}}\max\left(\frac{x}{p_{1}^{e_{1}}p_{2}^{e_{2}}\log(y)},1\right)<<\frac{x\log\log(y)}{\log(y)}.
$$
\noindent With regard to $\sum_{3}$, we have
$$
\sum_{3}=\sum_{p_{1},p_{2}<y^{1/4}}\Phi\left(\frac{x}{p_{1}p_{2}},y\right)+\sum_{\substack{p_{1},p_{2}<y\\
		\text{one of }p_{1},p_{2}>y^{1/4}}}\Phi\left(\frac{x}{p_{1}p_{2}},y\right)
$$
\begin{equation}
=\sum_{5}+\sum_{6}.
\end{equation}
\noindent Now, the latter sum in (17) is
\begin{equation}
\sum_{6}<<\sum_{\substack{p_{1},p_{2}<y\\
		\text{one of }p_{1},p_{2}>y^{1/4}}}\frac{x}{p_{1}p_{2}\log(y)}+\sum_{\substack{p_{1},p_{2}<y\\
		\text{one of }p_{1},p_{2}>y^{1/4}\\
		x/y<p_{1}p_{2}<x}}1,
\end{equation}
\noindent and the first sum in (18) may be bounded by
$$
\sum_{\substack{p_{1},p_{2}<y\\
		\text{one of }p_{1},p_{2}>y^{1/4}}}\frac{x}{p_{1}p_{2}\log(y)}<<\frac{x}{\log(y)}\left(\sum_{y^{1/4}<p<y}\frac{1}{p}\right)\left(\sum_{q<y}\frac{1}{q}\right)
$$
$$
<<\frac{x\log\log(y)}{\log(y)}.
$$
\noindent To bound the second sum in (18) note that if $p_{1}p_{2}>x/y=x^{1-1/\alpha}$, $\alpha>2$, then one of the primes $p_{1}, p_{2}>x^{1/4}$ and so
$$
\sum_{\substack{p_{1},p_{2}<y\\
		\text{one of }p_{1},p_{2}>y^{1/4}\\
		x/y<p_{1}p_{2}<x}}1<<\sum_{x^{1/4}<p_{1}<y}\pi\left(\frac{x}{p_{1}}\right)<<\frac{x}{\log(y)}\sum_{x^{1/4}<p_{1}<x}\frac{1}{p_{1}}<<\frac{x}{\log(y)}.
$$
\noindent Lastly,
$$
\sum_{5}=\sum_{p_{1},p_{2}<y^{1/t}}\frac{x}{\log(y)}\frac{w\left(\frac{\log(x)-\log(p_{1}p_{2})}{\log(y)}\right)}{p_{1}p_{2}}+O\left(\frac{x\log\log(y)}{\log^{2}(y)}\right)
$$
$$
=w(\alpha)\frac{x}{\log(y)}\sum_{p_{1},p_{2}<y^{1/t}}\frac{1}{p_{1}p_{2}}+O\left(\frac{x}{\log^{2}(y)}\sum_{p_{1},p_{2}<y^{1/t}}\frac{\log(p_{1}p_{2})}{p_{1}p_{2}}\right)
$$
\begin{equation}
=w(\alpha)\frac{x}{\log(y)}\frac{(\log\log(y^{1/4}))^{2}}{2}+O\left(\frac{x\log\log(y)}{\log(y)}\right)
\end{equation}
\noindent so that from (16), (17), (18), and (19), we have
$$
N_{2}(x,y)=w(\alpha)\frac{x(\log\log(y))^{2}}{2\log(y)}+O\left(\frac{x\log\log(y)}{\log(y)}\right).
$$
\indent This method can be used for $k>2$, the major difference in this case being that the primes in the above sum which were truncated at $y^{1/4}$ must now be truncated at $y^{1/t}$ with $t=2k$.  The method then yields

\bigskip

\noindent\textbf{Theorem 3*:} \textit{For fixed $k$ and fixed} $\alpha>2$
$$
N_{k}(x,y)=w(\alpha)\frac{x}{\log(y)}\frac{(\log\log(y))^{k}}{k!}+O\left(\frac{x(\log\log(y))^{k-1}}{(k-1)!\log(y)}\log(k)\right).
$$

\bigskip

\noindent\textbf{Remark:} In Theorem 3* the $\log(k)$ factor in the
error term is due to 
$$
\sum_{y^{1/2k}\leq p\leq y}\frac{1}{p} \asymp \log(k).$$

\indent We will see in Section 6 that by exploiting the methods of analysis we can derive an estimate of $N_{k}(x,y)$ for fixed $\alpha>2$ which is superior to Theorem 3.

\section{Estimate of $S_{z}(x,y)$ for Large $y$, for $\Re(z)>0$}

\indent This section will focus on estimates for $S_{z}(x,y)$ for \textit{large} values of $y$.  By \textit{large} we mean values of $y$ such that $y\geq\exp\left(\log^{1-\epsilon}(x)\right)$ for any sufficiently small $\epsilon>0$.  Our asymptotic estimates for these large values of $y$ are valid only when $\Re(z)>0$.  Nevertheless, a result from [13] (discussed later in the section) ensures that this is all the information that we need.

\indent Let us begin with the case $y \geq x$.  Clearly, in this case we have
\begin{equation}
S_{z}(x,y)=S_{z}(x,x)=\displaystyle\frac{g(1,z)}{\Gamma(z)}\frac{x}{\log^{1-z}(x)}+O\left(\frac{x}{\log^{2-z}(x)}\right)
\end{equation}
\noindent by the work of Selberg in [10].  We shall see that the limit
$$
m_{z}(\alpha):=\lim_{x\rightarrow\infty}\frac{S_{z}(x,y)}{x/\log^{1-z}(y)}
$$
\noindent exists for $\alpha>1$.  However, for $0<\alpha\leq1$ we may rewrite equation (20) as
$$
S_{z}(x,y)=m_{z}(\alpha)\frac{x}{\log^{1-z}(y)}+O_{\alpha}\left(\frac{x}{\log^{2-z}(y)}\right)
$$
\noindent where 
$$
m_{z}(\alpha)=\frac{g(1,z)}{\alpha^{1-z}\Gamma(z)}.
$$
\noindent With the aid of the following lemma we will be able to show that the function $m_{z}(\alpha)$ does in fact exist for larger ranges of $\alpha$.

\bigskip

\noindent\textbf{Lemma 5:} \textit{For $x\geq y^{h}\geq y\geq2$ and $|z|\leq R$ we have}
$$
S_{z}(x,y)=S_{z}(x,y^{h})+(1-z)\int_{y}^{x}S_{z}(x/t,t)\frac{dt}{\log(t)}+O(x R(y)\log^{R+1}(x)).
$$

\bigskip

\begin{proof} From the generalization of the sieve to strongly multiplicative
functions $g(n)$ as in [4],
we may choose $g(n)=z^{\omega_{y}(n)}$ to conclude that
\begin{equation}
S_{z}(x,y)=S_{z}(x,y^{h})+(1-z)\sum_{y\leq p\leq y^{h}}S_{z}(x/p,p).
\end{equation}
\noindent From the above equation we see that the lemma will be proven provided
$$
\sum_{y\leq p\leq y^{h}}S_{z}(x/p,p)=\int_{y}^{y^{h}}S_{z}(x/t,t)\frac{dt}{\log(t)}+O(xR(y)\log^{R+1}(x)).
$$
	
\noindent This is a standard exercise except that here 
we must decompose the intervals in a certain way.  When $p\leq n$ and $t\leq n$ we estimate the sum by decomposing the interval $[1,n]$ into the subintervals $I_{1}=[1,p_{1}],$ $I_{2}=(p_{1},p_{2}], ..., I_{\omega(n)}=(p_{\omega(n)-1},p_{\omega(n)}],$ $I_{\omega(n)+1}=(p_{\omega(n)}, n]$, where $p_{1}<p_{2}<...<p_{\omega(n)}$ denote the distinct prime divisors of $n$.  If $\chi_{j}(t)$ denotes the characteristic function of the interval $I_{j}$, then in this case the sum in (21) equals
$$
=\sum_{n\leq x/y}\sum_{j=1}^{\omega(n)+1}\left(\sum_{y<p<\min(x/n,y^{h})}\chi_{j}(n)z^{\omega_{p}(n)}-\int_{y}^{\min(x/n,y^{h})}\chi_{j}(t)z^{\omega_{t}(n)}\frac{dt}{\log(t)}\right)
$$
$$
=\sum_{n\leq x/y}\sum_{j=1}^{\omega(n)+1}z^{\omega_{p_{j}}(n)}\left(\sum_{y<p<\min(x/n,y^{h})}1-\int_{y}^{\min(x/n,y^{h})}\frac{dt}{\log(t)}\right)
$$
$$
<<\sum_{n\leq x/n}\sum_{j=1}^{\omega(n)+1}\tau_{R}(n)\frac{x}{n}R(y)=xR(y)\sum_{n\leq x/y}\frac{\tau_{R}(n)}{n}(\omega(n)+1)
$$
\begin{equation}
<<xR(y)\log(x)\sum_{n\leq x/y}\frac{\tau_{R}(n)}{n}<<xR(y)\log^{R+1}(x).
\end{equation}
	
\noindent The result then follows from (21) and (22).
	
\end{proof}

\indent We may now derive an estimate for the series $S_{z}(x,y)$ with
$1\leq\alpha\leq2$ using Buchstab's method.

\bigskip

\noindent\textbf{Theorem 4:} \textit{For $2\leq\sqrt{x}\leq y\leq x$, and $z\in\mathbb{C}$, $|z|\leq R$, $\Re(z)>0$, we have
$$
S_{z}(x,y)=m_{z}(\alpha)\frac{x}{\log^{1-z}(y)}+O\left(\frac{(1-z)x^{\ast}}{\log^{2-z}(y)}+\frac{(1-z)x^{\ast}}{\log(x)}\right),
$$
\noindent where $x^{\ast}:=x\log\log(x)$ and
\begin{equation}
m_{z}(\alpha)=\frac{g(1,z)}{\Gamma(z)}\left(\frac{1}{\alpha^{1-z}}+\frac{(1-z)}{\alpha^{1-z}}\int_{1}^{\alpha}\frac{du}{u^{z}(u-1)^{1-z}}\right)
\end{equation}
\noindent for $1\leq\alpha\leq2$}.

\bigskip

\begin{proof} We begin with (21) and take $y^{h}=x$ to obtain
\begin{equation}
S_{z}(x,y)=S_{z}(x,x)+(1-z)\sum_{y\leq p<x}S_{z}(x/p,p).
\end{equation}
\noindent As $S_{z}(x,x)=S_{z}(x)$ is already estimated by the Selberg asymptotic, we need only focus on the sum (which can be evaluated with the assistance of the previous lemma).  With this in mind we state the Selberg estimate for $S_{z}(x)$ in the slightly more convenient form of
\begin{equation}
S_{z}(x)=\frac{g(1,z)}{\Gamma(z)}\frac{x}{\max\left(\log(x),2\right)^{1-z}}+O\left(\frac{x}{\max\left(\log(x),2\right)^{2-z}}\right).
\end{equation}
\noindent Now, if $\sqrt{x}\leq y\leq x$ and $y\leq p<x$, then $x/p\leq\sqrt{x}\leq y\leq p$ so that the sum in equation (24) becomes
$$
\sum_{y\leq p<x}S_{z}(x/p,p)=\sum_{y\leq p<x}S_{z}(x/p)=\int_{y}^{x}S_{z}(x/t)\frac{dt}{\log(t)}+O(xR(y)\log^{R+1}(x))
$$
\noindent by the results in Lemma 6.  It therefore suffices to obtain an
accurate estimate of the above integral, which can be done by applying
Lemma 6 in the following form 
	
$$
\int_{y}^{x}S_{z}(x/t)\frac{dt}{\log(t)}=\frac{g(1,z)}{\Gamma(z)}\int_{y}^{x/e^{2}}\frac{x}{t\log^{1-z}(x/t)}\frac{dt}{\log(t)}
$$
\begin{equation}
+O\left(\int_{y}^{x/e^{2}}\frac{x}{t\log^{2-z}(x/t)}\frac{dt}{\log(t)}+\int_{x/e^{2}}^{x}\frac{x}{t}\frac{dt}{\log(t)}\right).
\end{equation}
\noindent Letting $t=x^{1/u}$ the first integral in equation (26) becomes
$$
\int_{y}^{x/e^{2}}\frac{x}{t\log^{1-z}(x/t)}\frac{dt}{\log(t)}=\int_{\left(1-\frac{2}{\log(x)}\right)^{-1}}^{\alpha}\frac{x\log(x)du}{u^{2}\left(\frac{\log(x)}{u}\right)\left(1-\frac{1}{u}\right)^{1-z}\log^{1-z}(x)}
$$
\begin{equation}	=\frac{x}{\log^{1-z}(x)}\int_{\left(1-\frac{2}{\log(x)}\right)^{-1}}^{\alpha}\frac{du}{u^{z}(u-1)^{1-z}}.
\end{equation}
\noindent  Also, the second integral in the $O$-term in (26) is 

$$
\int_{x/e^{2}}^{x}\frac{x}{t}\frac{dt}{\log(t)}=x\int_{x/e^{2}}^{x}\frac{dt}{t\log(t)}
$$
$$
=x\left(\log\log(x)-\log\log(x/e^{2})\right)
$$
\begin{equation}
=-x\log\left(1-\frac{2}{\log(x)}\right)<<\frac{x}{\log(x)},.
\end{equation}
\noindent  We see from equations (26), (27), and (28) that
$$
\int_{y}^{x}S_{z}(x/t)\frac{dt}{\log(t)}=
$$
\begin{equation}
\frac{g(1,z)}{\Gamma(z)}\frac{x}{\log^{1-z}(x)}\int_{\left(1-\frac{2}{\log(x)}\right)^{-1}}^{\alpha}\frac{du}{u^{z}(u-1)^{1-z}}+O\left(I_{2}+\frac{x}{\log(x)}\right).
\end{equation}
\noindent where
$$
I_{2}:=\int_{y}^{x/e^{2}}\frac{x}{t\log^{2-z}(x/t)}\frac{dt}{\log(t)}.
$$
\noindent We will now estimate the integral in the main term of equation (29), as
\begin{equation}
\int_{\left(1-\frac{2}{\log(x)}\right)^{-1}}^{\alpha}\frac{du}{u^{z}(u-1)^{1-z}}=\int_{1}^{\alpha}\frac{du}{u^{z}(u-1)^{1-z}}-\int_{1}^{\left(1-\frac{2}{\log(x)}\right)^{-1}}\frac{du}{u^{z}(u-1)^{1-z}}.
\end{equation}
\noindent Note that $1\leq\alpha\leq2$ so that if $\Re(z)>0$, then 
$$
\left|\int_{1}^{\alpha}\frac{(u-1)^{z-1}}{u^{z}}du\right|\leq\int_{1}^{\alpha}\frac{du}{u^{\Re(z)}}\int_{1}^{\alpha}(u-1)^{\Re(z)-1}du
$$
$$
=\frac{(\alpha-1)^{\Re(z)}}{\Re(z)}\int_{1}^{\alpha}\frac{du}{u^{\Re(z)}}<+\infty,
$$
\noindent that is, the integrals in (30) are convergent, provided $\Re(z)>0$
(here is where we must use the fact that the real part of $z$ is positive,
as the integral in (30) would not converge if $\Re(z)\leq0$).  We wish to
bound the second integral in equation (30). To that end
$$
\int_{1}^{\left(1-\frac{2}{\log(x)}\right)^{-1}}\frac{du}{u^{z}(u-1)^{1-z}}<<\int_{1}^{\left(1-\frac{2}{\log(x)}\right)^{-1}}\frac{du}{(u-1)^{1-z}}<<\frac{1}{\log^{\Re(z)}(x)};
$$
\noindent and consequently
\begin{equation}
\int_{\left(1-\frac{2}{\log(x)}\right)^{-1}}^{\alpha}\frac{du}{u^{z}(u-1)^{1-z}}=\int_{1}^{\alpha}\frac{du}{u^{z}(u-1)^{1-z}}+O\left(\frac{1}{\log^{\Re(z)}(x)}\right).
\end{equation}
	
\indent The evaluation of the integral $I_{2}$ is similar to what we have done for the main term, except here we must show care when $\Re(z)$ is close to $1$.  This is achieved by considering the two cases $|\Re(z)-1|>\frac{1}{5\log\log(x)}$ and $|\Re(z)-1|\leq\frac{1}{5\log\log(x)}$.  The substitution $t=x^{1/u}$ shows that
$$
I_{2}=O\left(\frac{x\log\log(x)}{\log^{2-z}(x)}\right)=O\left(\frac{x^{\ast}}{\log^{2-z}(x)}\right)
$$
\noindent for $\Re(z)\geq1$ and
$$
I_{2}=O\left(\frac{x\log\log(x)}{\log(x)}\right)=O\left(\frac{x^{\ast}}{\log(x)}\right)
$$
\noindent when $\Re(z)<1$, and where we have set $x^{\ast}=x\log\log(x)$; thus
$$
\int_{y}^{x}S_{z}(x/t)\frac{dt}{\log(t)}
$$
$$
=\frac{g(1,z)}{\Gamma(z)}\frac{x}{\log^{1-z}(x)}\int_{1}^{\alpha}\frac{du}{u^{z}(u-1)^{1-z}}+O\left(\frac{x^{\ast}}{\log(x)}+\frac{x^{\ast}}{\log^{2-z}(x)}\right),
$$
\noindent from (31) and these estimates the result follows.
	
\end{proof}

\noindent  We may also note that if $r\in\mathbb{R}^{+}$ then the preceding result can be stated without the $\log\log(x)$ factors in $x^{\ast}$, that is

\bigskip

\noindent\textbf{Theorem 4*:} \textit{For $2\leq\sqrt{x}\leq y\leq x$, and $r\in\mathbb{R}^{+}$, $r\leq R$, we have
$$
S_{r}(x,y)=m_{r}(\alpha)\frac{x}{\log^{1-r}(y)}+O\left(\frac{x}{\log^{2-r}(y)}+\frac{x}{\log(x)}\right),
$$
\noindent where $m_{r}(\alpha)$ is given by (23).}

\bigskip

\indent We may now use Buchstab's recurrence (equation (24)) to derive an estimate for $S_{z}(x,y)$ for $\alpha\geq2$ by applying induction on $[\alpha]$.  Fortunately, we will also deduce an asymptotic estimate for $S_{z}(x,y)$ when $\Re(z)>0$ for $y$ outside the range of Corollary 1.  First we prove

\bigskip

\noindent\textbf{Theorem 5:} \textit{Let $\Re(z)>0$. Then for arbitrary but fixed $\alpha>2$, we have
$$
S_z(x,y)=\frac{xm_{z}(\alpha)}{log^{1-z}y} +O\left(\frac{x^{\ast}}{log^{2-z}(y)}+\frac{x^{\ast}}{log(y)}\right),
$$
\noindent where for $\alpha>2$, $m_{z}(\alpha)$ is given by}
\begin{equation}
m_{z}(\alpha)=\frac{2^{1-z}m_{z}(2)}{\alpha^{1-z}}+\frac{1-z}{\alpha^{1-z}}\int_{2}^{\alpha}\frac{m_{z}(u-1)}{u^{z}}du.
\end{equation}

\bigskip

\begin{proof} Setting $y^{h}=\sqrt{x}$ in equation (21)
$$
S_{z}(x,y)=S_{z}(x,\sqrt{x})+(1-z)\sum_{y\leq p\leq\sqrt{x}}S_{z}(x/p,p)
$$
$$
=\frac{2^{1-z}m_{z}(2)}{\alpha^{1-z}}\frac{x}{\log^{1-z}(y)}+O\left(\frac{(1-z)x^{\ast}}{\log^{2-z}(y)}\right)
$$
$$
+(1-z)\int_{y}^{\sqrt{x}}m_{z}\left(\frac{\log(x/t)}{\log(t)}\right)\frac{x/t}{\log^{1-z}(t)}\frac{dt}{\log(t)}
$$
$$+O\left((1-z)\int_{y}^{\sqrt{x}}\frac{x^{\ast}/t}{\log^{2-z}(t)}\frac{dt}{\log(t)}\right)
+O\left((1-z)\int_{y}^{\sqrt{x}}\frac{x^{\ast}}{t\log^{2}(t)}dt\right)
$$
$$
+O(xR(y)\log^{R}(x)),
$$
\noindent as $\sqrt{x}=y^{\alpha/2}$.  With the familiar substitution of $t=x^{1/u}$ we obtain
$$
=\frac{2^{1-z}m_{z}(2)}{\alpha^{1-z}}\frac{x}{\log^{1-z}(y)}+O_{R,\alpha}\left(\frac{x^{\ast}}{\log^{2-z}(y)}\right)
$$
$$
+\frac{(1-z)}{\alpha^{1-z}}\frac{x}{\log^{1-z}(y)}\int_{2}^{\alpha}m_{z}(u-1)\frac{du}{u^{z}}+O_{R,\alpha}\left(\frac{x^{\ast}}{\log^{2-z}(y)}\right)+O\left(\frac{x^{\ast}}{\log(y)}\right);
$$
\noindent hence, the theorem follows from this with $m_{z}(\alpha)$ as defined in (32) above.

\end{proof}

\indent From the definition of $m_{z}(\alpha)$ in (32), we see that
$$
\alpha^{1-z}m_{z}(\alpha)-(\alpha-1)^{1-z}m_{z}(\alpha-1)=(1-z)\int_{\alpha-1}^{\alpha}\frac{m_{z}(u-1)}{u^{z}}du,
$$
\noindent for $\alpha\geq3$.  This can be rewritten as:
$$
m_{z}(\alpha):=\frac{m_{z}(\alpha-1)(\alpha-1)^{1-z}}{\alpha^{1-z}}+\frac{(1-z)}{\alpha^{1-z}}\int_{\alpha-1}^{\alpha}m_{z}(u-1)\frac{du}{u^{z}}.
$$

\indent We next derive an improvement of Theorem 4 in which the asymptotic estimate will hold for $\alpha$ tending to infinity with $x$.  For this (32) will be useful.

\bigskip

\noindent\textbf{Theorem 6:} \textit{Let $\alpha\geq1$, if $\Re(z)\geq1$ there exists an absolute constant $K=K(R)$ such that
$$
\left|S_{z}(x,y)-m_{z}(\alpha)\frac{x}{\log^{1-z}(y)}\right|<<\frac{\alpha^{K}x^{\ast}}{\log^{2-z}(y)},
$$
\noindent and if $0<\Re(z)<1$,  then}
$$
\left|S_{z}(x,y)-m_{z}(\alpha)\frac{x}{\log^{1-z}(y)}\right|\leq\frac{\alpha^{K}x^{\ast}}{\log(y)}.
$$

\bigskip

\begin{proof} Clearly, it suffices to prove the theorem for $\alpha>3$.  For $\alpha>3$, we can use (32).  Let $y=x^{1/\alpha}$, $y^{h}=x^{1/(\alpha-1)}$, and $u=\frac{\log(x)}{\log(t)}$.  Lemma 6 shows that, with this notation,
$$
S_{z}(x,y)=
$$
\begin{equation}
S_{z}(x,y^{h})+(1-z)\int_{y}^{y^{h}}S_{z}(x/t,t)\frac{dt}{\log(t)}+O\left(x\alpha^{|z|+1}R(y)\log^{R+1}(y)\right).
\end{equation}
\noindent  Assume that $\Re(z)\geq1$.  We shall prove the result by induction on $\alpha$.  Assume that there exists a positive, non-decreasing function $\phi(u)$ such that for all $u\leq\alpha-1$ and $x>y>1$ we have
\begin{equation}
\left|S_{z}(x,t)-m_{z}(u)\frac{x}{\log^{1-z}(t)}\right|<\frac{\phi(u)x^{\ast}}{\log^{2-z}(t)},
\end{equation}
\noindent as $\Re(z)\geq1$ Theorem 4 establishes the validity of (34) when $\alpha\in[2,3]$.  By the inductive hypothesis of (34) and equation (33) we now obtain
$$
S_{z}(x,y)=m_{z}(\alpha-1)\left(\frac{\alpha-1}{\alpha}\right)^{1-z}\frac{x}{\log^{1-z}(y)}+O_{1}\left(\left(\frac{\alpha-1}{\alpha}\right)^{2-z}\frac{\phi(\alpha-1)x^{\ast}}{\log^{2-z}(y)}\right)
$$
$$
+(1-z)\int_{y}^{y^{h}}m_{z}\left(\frac{\log(x)-\log(t)}{\log(t)}\right)\frac{x}{t\log^{2-z}(t)}dt
$$
$$
+O_{1}\left((1-z)\int_{y}^{y^{h}}\phi\left(\frac{\log(x)-\log(t)}{\log(t)}\right)\frac{x^{\ast}}{t\log^{3-z}(t)}\right)dt
$$
\begin{equation}
+O(xR(y)\alpha^{|z|}\log^{R+1}(y)),
\end{equation}
	
\noindent where the notation $O_{1}$ implies that the implicit constant is $\leq1$ (which will be important when iterating the process).  We will see that the first and third terms in the above recurrence will make the largest contribution, and note that as $t=x^{1/u}$
$$
m_{z}(\alpha-1)\left(\frac{\alpha-1}{\alpha}\right)^{1-z}\frac{x}{\log^{1-z}(y)}+(1-z)\int_{y}^{y^{h}}m_{z}\left(\frac{\log(x)-\log(t)}{\log(t)}\right)\frac{x}{t\log^{2-z}(t)}
$$
$$
=m_{z}(\alpha-1)\left(\frac{\alpha-1}{\alpha}\right)^{1-z}\frac{x}{\log^{1-z}(y)}+\frac{(1-z)}{\alpha^{1-z}}\frac{x}{\log^{1-z}(y)}\int_{\alpha-1}^{\alpha}m_{z}\left(u-1\right)\frac{du}{u^{z}}
$$

\begin{equation}
=m_{z}(\alpha)\frac{x}{\log^{1-z}(y)},
\end{equation}
\noindent from (32).  We conclude from (35) and (36) that
$$
S_{z}(x,y)=\frac{m_{z}(\alpha)x}{\log^{1-z}(y)}+O_{1}\left(\left(\frac{\alpha-1}{\alpha}\right)^{2-z}\frac{\phi(\alpha-1)x^{\ast}}{\log^{2-z}(y)}\right)
$$
\begin{equation}
+O_{1}\left((1-z)x^{\ast}\phi(\alpha-1)\int_{y}^{y^{h}}\frac{dt}{t\log^{3-z}(t)}\right)+O(x\alpha^{|z|+1}R(y)\log^{R+1}(y)).
\end{equation}
\noindent If $\Re(z)\neq2$ the integral in equation (37) will be equal to
$$
\int_{y}^{y^{h}}\frac{dt}{t\log^{3-\Re(z)}(t)}=\frac{\log^{\Re(z)-2}(y^{h})}{\Re(z)-2}-\frac{\log^{\Re(z)-2}(y)}{\Re(z)-2}
$$
$$
=\frac{1}{(\Re(z)-2)}\frac{1}{\log^{2-\Re(z)}(y)}\left(1-\frac{1}{h^{2-\Re(z)}}\right),
$$
\noindent and
\begin{equation}
1-\frac{1}{h^{2-\Re(z)}}=1-\left(1-\frac{1}{\alpha}\right)^{2-\Re(z)}=O\left(\frac{|2-\Re(z)|}{\alpha}\right).
\end{equation}
\noindent If $\Re(z)=2$, then the integral is
$$
\int_{y}^{y^{h}}\frac{dt}{t\log^{3-\Re(z)}(t)}=\int_{y}^{y^{h}}\frac{dt}{t\log(t)}=\log\log(y^{h})-\log\log(y)
$$
$$
=\log(h)+\log\log(y)-\log\log(y)=\log(h)
$$
\begin{equation}
=\log\left(\frac{\alpha}{\alpha-1}\right)=1+O\left(\frac{1}{\alpha}\right).
\end{equation}
\noindent From equations (37), (38), and (39) equation (35) becomes
$$
S_{z}(x,y)=\frac{m_{z}(\alpha)x}{\log^{1-z}(y)}+O_{1}\left(\left(\frac{\alpha-1}{\alpha}\right)^{2-z}\frac{\phi(\alpha-1)x^{\ast}}{\log^{2-z}(y)}\right)
$$
\begin{equation}
+O_{1}\left(\frac{|1-z|x^{\ast}\phi(\alpha-1)}{\log^{2-z}(y)}\frac{C}{\alpha}\right).
\end{equation}
\noindent We are free to choose our function $\phi(u)$ to satisfy
$$
\phi(\alpha)<\phi(\alpha-1)\left(\frac{\alpha-1}{\alpha}\right)^{1-z}+\phi(\alpha-1)\frac{|1-z|C}{\alpha}\leq\phi(\alpha-1)\left(1+\frac{|1-z|C^{\prime}}{\alpha}\right).
$$
	
\indent Because the function $\phi(\alpha)$ is defined recursively we may estimate its growth rate by noting that for any $j\in\mathbb{Z}$, $0\leq j\leq\alpha-2$
$$
\frac{\phi(\alpha-j)}{\phi(\alpha-j-1)}<1+\frac{|1-z|C}{\alpha-j};
$$
\noindent hence,
$$
\phi(\alpha)=\frac{\phi(\alpha)}{\phi(\alpha-1)}\frac{\phi(\alpha-1)}{\phi(\alpha-2)}...\frac{\phi(2)}{\phi(1)}<\prod_{0\leq j\leq\alpha-1}\left(1+\frac{|1-z|C}{\alpha-j}\right)
$$
$$
=\prod_{1\leq j\leq \alpha}\left(1+\frac{|1-z|C}{j}\right)<<\alpha^{C|1-z|},
$$
\noindent from the estimate
$$
\prod_{1\leq j\leq x}\left(1+\frac{C}{j}\right)<<x^{C}.
$$
\noindent Collecting the above results we may conclude that
$$
\left|S_{z}(x,y)-m_{z}(\alpha)\frac{x}{\log^{1-z}(y)}\right|<\frac{\alpha^{C|1-z|}x^{\ast}}{\log^{2-z}(y)},
$$
\noindent for some constant $K=K(R)=C|1-z|$, thereby proving the first statement of Theorem 6.
	
\indent If $0<\Re(z)<1$ then the method of Alladi in [1] shows that
$$
\left|S_{z}(x,y)-m_{z}(\alpha)\frac{x}{\log^{1-z}(y)}\right|\leq\frac{\alpha x^{\ast}}{\log(y)}
$$
\noindent in place of (34).  The second statement of the theorem then follows by repeating the above induction procedure.
	
\end{proof}

\indent The preceding theorem gives the desired uniform result we need to analyze the sum $S_{z}(x,y)$ for large values of $y$.  Let $\Re(z)\geq\delta>0$ and consider

$$
\log^{\delta}(y)=\alpha^{K}\log\log(x)=\left(\frac{\log(x)}{\log(y)}\right)^{K}\log\log(x)\leq\frac{\log^{K+\epsilon}(x)}{\log^{K}(y)}
$$
\noindent for some $0<\epsilon<\delta$ so that
$$
\log(y)=\log^{(k+\epsilon)/(k+\delta)}(x)=\log^{1-\delta^{\prime}}(x)
$$
\noindent and
$$
y=\exp(\log^{1-\delta^{\prime}}(x)).
$$
\noindent Therefore, we obtain an asymptotic estimate of $S_{z}(x,y)$ from Theorem 6 when $y\geq e^{\log^{1-\delta^{\prime}}(x)}$, which is larger than the range specified in Corollary 1.

\indent Let us pause for a moment and reflect on what has just been proven.  In Theorem 1 we obtained an asymptotic formula for $S_{z}(x,y)$ ($|z|\leq R$), provided $y\leq x^{1/(R+D+1+\epsilon)\log\log(x)}$; however, Theorem 6 gives an asymptotic formula (for $0<\Re(z)$, $|z|\leq R$) provided $y\geq x^{1/\log^{\delta}(x)}$.
Observe that for $x$ sufficiently large, we have
$$
x^{1/\log^{\delta}(x)}\leq x^{1/(R+D+1+\epsilon)\log\log(x)}
$$
\noindent so the ranges of Corollary 2 and Theorem 6 overlap.  We have by
virtue of Theorem 6 and Corollary 2 derived an asymptotic formula for all
$y\leq x$ if $R\geq\Re(z)>0$ (note, however, that the results of Corollary 2
are true with only the restriction $|z|\leq R$).  This uniform estimate
will be one of the main tools utilized in Section 6 for the study of the
local distribution of small prime factors.

\section{Properties of $m_{z}(\alpha)$}

\indent In this section we shall study the properties of the function
$m_{z}(\alpha)$ arising in Theorem 5.  We have already shown that this
function exists and is given by (32). Thus

$$
\alpha^{1-z}m_{z}(\alpha)=2^{1-z}m_{z}(2)+(1-z)\int_{2}^{\alpha}m_{z}(u-1)\frac{du}{u^{z}};
$$
\noindent so that $m_{z}(\alpha)$ satisfies the following
difference-differential equation:
\begin{equation}
(\alpha^{1-z}m_{z}(\alpha))^{\prime}=(1-z)\left(\frac{m_{z}(\alpha-1)}{\alpha^{z}}\right),
\end{equation}
\noindent where $m_{z}^{\prime}(\alpha)=\frac{d}{d\alpha}m_{z}(\alpha)$. Equation
(41) implies that
$$
\alpha^{1-z}m_{z}^{\prime}(\alpha)+(1-z)\alpha^{-z}m_{z}(\alpha)=(1-z)m_{z}(\alpha-1)\alpha^{-z}
$$
\noindent which can be rewritten as
\begin{equation}
m_{z}^{\prime}(\alpha)=\frac{1-z}{\alpha}\left(m_{z}(\alpha-1)-m_{z}(\alpha)\right)=\frac{z-1}{\alpha}\int_{\alpha-1}^{\alpha}m_{z}^{\prime}(u)du.
\end{equation}

\noindent The above equation will allow us to derive some useful
properties about the convergence of $m_{z}(\alpha)$ as
$\alpha\rightarrow\infty$.

\bigskip

\noindent\textbf{Lemma 6:} \textit{For $z\in\mathbb{C}$, $\Re(z)>0$, $|z|\leq R$, and $\alpha\geq2$ the function $m_{z}^{\prime}(\alpha)$ is differentiable.
Moreover,} $m_{z}^{\prime}(\alpha)<<e^{-\alpha\log\alpha+O_{R}(\alpha)}$.

\bigskip

\begin{proof} We may immediately deduce from equation (42) that $m_{z}^{\prime}(\alpha)$ is differentiable.  We will first show that equation (42) implies that $m_{z}^{\prime}(\alpha)$ is uniformly bounded for $\alpha\in\mathbb{R}$, $\alpha\geq2$.
	
\indent Suppose not.  Then $|m_{z}^{\prime}(\alpha)|$ assumes all sufficiently large values by continuity.  Define
$$
B:=\max_{\alpha\leq|z-1|+1}|m_{z}^{\prime}(\alpha)|
$$
\noindent and consider a value $M>B$ assumed by $|m_{z}^{\prime}(\alpha)|$.
Next, define
$$
\alpha_{0}:=\inf\{\alpha: |m_{z}(\alpha)|=M\}
$$
\noindent so that $|m_{z}^{\prime}(\alpha_{0})|=M$ and $\alpha_{0}>|z-1|+1$.  Now, from (42) we see
$$
M=|m_{z}^{\prime}(\alpha_{0})|\leq\frac{|z-1|}{\alpha_{0}}\int_{\alpha_{0}-1}^{\alpha_{0}}|m_{z}^{\prime}(t)|dt<\frac{|z-1|}{\alpha_{0}}M<M,
$$
\noindent a contradiction.  It follows that	
$$
\{\alpha: |m_{z}(\alpha)|=M\}=\varnothing
$$
\noindent for each $M>B$, so that $\max|m_{z}^{\prime}(\alpha)|\leq B$.
	
\indent Using the above analysis, and again invoking equation (42), we see that
$$
|m_{z}^{\prime}(\alpha)|=\frac{|1-z|}{\alpha}\left|\int_{\alpha-1}^{\alpha}m_{z}^{\prime}(u)du\right|\leq\frac{|1-z|}{\alpha}\int_{\alpha-1}^{\alpha}|m_{z}^{\prime}(u)|du
<\frac{B|1-z|}{\alpha}.
$$
\noindent Consider the function $\sup_{u\in[t,+\infty)}|m_{z}^{\prime}(u)|<N$ which is clearly monotone and non-increasing, so
$$
|m_{z}^{\prime}(\alpha)|\leq\frac{|1-z|}{\alpha}\int_{\alpha-1}^{\alpha}\sup_{u\in[t,+\infty)}|m_{z}^{\prime}(u)|dt\leq\frac{|1-z|}{\alpha}\sup_{t\in[\alpha-1,+\infty)}|m_{z}^{\prime}(t)|
$$
$$
<\frac{|1-z|}{\alpha}\left(\frac{|1-z|}{\alpha-1}B\right)=\frac{|1-z|^{2}}{\alpha(\alpha-1)}B
$$
\noindent and by iteration
\begin{equation}
|m_{z}^{\prime}(\alpha)|\leq\frac{|1-z|^{[\alpha]}}{\Gamma([\alpha]+1)}B.
\end{equation}
\noindent Since $\frac{1}{\Gamma(\alpha+1)}<<e^{-\alpha\log\alpha}$, we deduce
from (43) that
$$
m_{z}^{\prime}(\alpha)<<e^{-\alpha\log\alpha}(1+R)^{[\alpha]}=e^{-\alpha\log\alpha}e^{[\alpha]\log(1+R)}=e^{-\alpha\log\alpha+O_{R}(\alpha)}.
$$

\end{proof}

\indent The rapid rate of decay of $m_{z}^{\prime}(\alpha)$ as
$\alpha\rightarrow\infty$ forces the function $m_{z}(\alpha)$ to approach a
limit $\ell(z)$.  We have already observed in Section 2 that for certain
ranges of $\alpha$, the estimate for the sum $S_{z}(x,y)$ is supplied by
Theorem 1.  The fact that the ranges of our estimates overlap will allow us
to obtain a representation for the limiting function $\ell(z)$.

\bigskip

\noindent\textbf{Theorem 7:} \textit{There exists a function $\ell(z)$ such that
$$
m_{z}(\alpha)=\ell(z)+O(e^{-\alpha\log\alpha+O_{R}(\alpha)})
$$
\noindent as $\alpha\rightarrow\infty$, and is given by
$$
\ell(z)=e^{(z-1)\gamma}\prod_{p}\left(1+\frac{z-1}{p}\right)\left(1-\frac{1}{p}\right)^{z-1}.
$$
\noindent Thus $\ell(z)$ given by the product is an analytic function of
$z\in\mathbb{C}$, and is nonzero if $z\neq1$ or $1-p$, where $p$ denotes
a prime number.}

\bigskip

\begin{proof}
We have seen in Lemma 7 that 
\begin{equation}
\frac{d}{d\alpha}m_{z}(\alpha)=m_{z}^{\prime}(\alpha)<<\frac{|z-1|^{\alpha}}{\Gamma(\alpha)}
\end{equation}
\noindent so that
\begin{equation}
f(z):=\int_{1}^{\infty}m_{z}^{\prime}(t)dt
\end{equation}
\noindent exists, as equation (44) implies the convergence of the integral in (45).  Now,
$$
\int_{1}^{\alpha}m_{z}^{\prime}(t)dt=m_{z}(\alpha)-m_{z}(1)
$$
\noindent so that
\begin{equation}
\int_{\alpha}^{\infty}m_{z}^{\prime}(t)dt=\int_{1}^{\infty}m_{z}^{\prime}(t)dt-\int_{1}^{\alpha}m_{z}^{\prime}(t)dt=f(z)-m_{z}(\alpha)+m_{z}(1).
\end{equation}
\noindent However,
\begin{equation}
\left|\int_{\alpha}^{\infty}m_{z}^{\prime}(t)dt\right|<<\int_{\alpha}^{\infty}e^{-t\log(t)+O(t)}dt<<_{R}e^{-\alpha\log\alpha+O_{R}(\alpha)}
\end{equation}
	
\noindent so that from (46) and (47) we obtain
$$
f(z)-m_{z}(\alpha)+m_{z}(1)<<_{R}e^{-\alpha\log\alpha+O_{R}(\alpha)}
$$
\noindent or
$$
m_{z}(\alpha)=f(z)+m_{z}(1)+O_{R}(e^{-\alpha\log\alpha+O_{R}(\alpha)})
$$
\noindent thereby proving the limit exists with
$$
\ell(z)=f(z)+m_{z}(1)=\int_{1}^{\infty}m_{z}^{\prime}(t)dt+m_{z}(1).
$$
	
\indent Now, recall from the comments following the proof of Theorem 6 that if we choose $(R+D+2)\log\log(x)<\alpha<(R+D+3)\log\log(x)$ then we may apply both Theorem 6 and Corollary 1; therefore, we need only equate $\ell(z)$ with the corresponding term in Corollary 1.  For $\alpha\asymp\log\log(x)$ Corollary 1 supplies
$$
\frac{S_{z}(x,y)}{x/\log^{1-z}(y)}=\log^{1-z}(y)\prod_{p<y}\left(1+\frac{z-1}{p}\right)+O\left(\frac{1}{\log^{\epsilon}(y)}\right)
$$
$$
=e^{(z-1)\gamma}\prod_{p<y}\left(1+\frac{z-1}{p}\right)\left(1-\frac{1}{p}\right)^{z-1}+O\left(\frac{1}{\log^{\epsilon}(y)}\right),
$$

\noindent so that as $\alpha\rightarrow\infty$ we see that
$$
\ell(z)=e^{(z-1)\gamma}\prod_{p}\left(1+\frac{z-1}{p}\right)\left(1-\frac{1}{p}\right)^{z-1}.
$$	
	
\end{proof}

\indent Having shown that the limit $\lim_{\alpha\rightarrow\infty}m_{z}(\alpha)=\ell(z)$ exists, it is clear from the product formula of $\ell(z)$ that the limit is nonzero for $z\neq1$ or $1-p$.  However, for future applications we must demonstrate that $m_{r}(\alpha)>0$ for all $r>0$.  With this in mind, we now prove the following corollary.

\bigskip

\noindent\textbf{Corollary 3:} \textit{The function $m_{r}(\alpha)>0$ for all real $r\geq0$ and $\alpha\geq1$.}

\bigskip

\begin{proof} Let $0\leq r<1$ so that we may immediately conclude from (32) that $m_{r}(1)>0$.  Now we may proceed by induction on $[\alpha]$ and utilize (32) to see that $m_{r}(\alpha>0)$ for all $\alpha\geq1$.
	
\indent If $r=1$ then the function $m_{1}(\alpha)\equiv1$, so we need only consider the case when $r$ is greater than $1$ for the corollary to be proven.  However, if $r>1$ then $S_{r}(x,y)=\sum_{n\leq x}r^{\omega_{y}(n)}$ will be an increasing function of $y$ and, consequently $S_{r}(x,y)$ will be a decreasing function of $\alpha$.  Let us define the function $m_{r}^{\ast}(\alpha)$ to satisfy
$$
S_{r}(x,y)=m_{z}^{\ast}(\alpha)\frac{x}{\log^{1-r}(x)}+O\left(\frac{x^{\ast}}{\log^{2-r}(x)}+\frac{x^{\ast}}{\log(x)}\right)
$$
\noindent and so $m_{r}^{\ast}(\alpha)\alpha^{r-1}=m_{r}(\alpha)$, with $m_{r}^{\ast}(\alpha)$ decreasing.  As $S_{r}(x,y)>0$ we see that $m_{r}^{\ast}(\alpha)>0$ for all $\alpha\geq1$, and this forces $m_{r}(\alpha)>0$ for all $\alpha\geq1$.

\end{proof}

\bigskip

\noindent\textbf{Theorem 8:} \textit{For each $\alpha>1$}
$$
\lim_{r\rightarrow0^{+}}m_{r}(\alpha)=w(\alpha).
$$

\bigskip

\begin{proof} Note that for $r>0$, $1<\alpha<2$, $m_{r}(\alpha)$ is given by 
$$
m_{r}(\alpha)=\frac{s(r)}{\alpha^{1-r}}+\frac{(1-r)s(r)}{\alpha^{1-r}}\int_{1}^{\alpha}\frac{du}{u^{r}(u-1)^{1-r}}
$$
\noindent (where $s(r)=g(1,r)/\Gamma(r)$), and for $\alpha>2$
\begin{equation}
\left(\alpha^{1-r}m_{r}(\alpha)\right)^{\prime}=(1-r)\frac{m_{r}(\alpha-1)}{\alpha^{r}}.
\end{equation}
\noindent Clearly when $r\rightarrow0$, the differential equation in (48)
coincides with that of the Buchstab function given in Section 3.  So we only
prove that $\lim_{r\rightarrow0^{+}}m_{r}(\alpha)=w(\alpha)$, for $1<\alpha<2$.
So we must show that
$$
\lim_{r\rightarrow0^{+}}\left(\frac{s(r)}{\alpha^{1-r}}+\frac{(1-r)s(r)}{\alpha^{1-r}}\int_{1}^{\alpha}\frac{du}{u^{r}(u-1)^{1-r}}\right)=\frac{1}{\alpha}.
$$
\noindent Since $\lim_{r\rightarrow0^{+}}s(r)=0$, we must show that
\begin{equation}
\lim_{r\rightarrow0^{+}}\frac{(1-r)s(r)}{\alpha^{1-r}}\int_{1}^{\alpha}\frac{du}{u^{r}(u-1)^{1-r}}=\frac{1}{\alpha}.
\end{equation}
\noindent Note that $s(r)\sim r$ as $r\rightarrow0^{+}$, because $g(1)=1$ and $\Gamma(r)$ has a simple pole with residue $1$ at $r=0$.  Using $u^{-r}=e^{-r\log(u)}=1+O(r)$, for $1<u<\alpha<2$ and $r\rightarrow0^{+}$, we have
$$
\int_{1}^{\alpha}\frac{du}{u^{r}(u-1)^{1-r}}=\int_{1}^{\alpha}\left(1+O(r)\right)\frac{du}{(u-1)^{1-r}}
$$
\begin{equation}
=\int_{1}^{\alpha}\frac{du}{(u-1)^{1-r}}+O\left(r\int_{1}^{\alpha}\frac{du}{(u-1)^{1-r}}\right)=\frac{(\alpha-1)^{r}}{r}+O((\alpha-1)^{r}).
\end{equation}
\noindent From (50) we obtain
\begin{equation}
\frac{(1-r)s(r)}{\alpha^{1-r}}\int_{1}^{\alpha}\frac{du}{u^{r}(u-1)^{1-r}}=\frac{(1-r)s(r)}{\alpha^{1-r}}\left(\frac{(\alpha-1)^{r}}{r}+O((\alpha-1)^{r})\right),
\end{equation}
\noindent and taking the limit as $r\rightarrow0^{+}$ in (51) gives the limit
in (49) as 
$$
=\lim_{r\rightarrow0^{+}}\frac{(1-r)s(r)}{\alpha^{1-r}}\left(\frac{(\alpha-1)^{r}}{r}+O((\alpha-1)^{r})\right)=\frac{1}{\alpha}\lim_{r\rightarrow0^{+}}\frac{s(r)}{r}=\frac{1}{\alpha}.
$$
\noindent This proves the theorem.

\end{proof}

{\textbf{Remark:}} We had defined $m_z(\alpha)$ only for $Re(z)>0$, but
$S_0(x,y)=\Phi(x,y)$. Thus it is important to establish Theorem 8, from which
we could interpret $w(\alpha)$ as $m_0(\alpha)$.

\section{The Local Distribution of the Number of Small Prime Factors}

\indent In this section we will apply the analytic results obtained in section's 2, 4, and 5 to study the function $N_{k}(x,y)$ by the contour integral method of Selberg in [10], by which we mean that we will apply the Cauchy integral formula to $S_{z}(x,y)$ in the following form
\begin{equation}
N_{k}(x,y)=\frac{1}{2\pi i}\int_{\partial U_{R}}\frac{S_{z}(x,y)}{z^{k+1}}dz,
\end{equation}
\noindent $U_{r}$ being the circle of radius $r$ centered at the origin.  We will begin with an analysis of the results which can be obtained from Theorem 1, which holds uniformly for \textit{small} $y$ and $|z|\leq R$.  We will then study $N_{k}(x,y)$ for \textit{large} values of $y$, in which case our results only apply for $\Re(z)>0$.  However, as was alluded to at the end of section 4, Theorem 11 due to Tenenbaum [13] will allow us to estimate $N_{k}(x,y)$ for \textit{large} values of $y$ provided we have an estimate for $S_{r}(x,y)$ for real $r>0$, which is supplied by the results of section 4.

\indent Tenenbaum [14] has also supplied the following alternative
representation for the function $m_{r}(\alpha)$:

\bigskip

\noindent\textbf{Theorem 9:} \textit{We have
$$
m_{r}(\alpha)=C(r)\left(\int_{0}^{\alpha-1}w(\alpha-t)\rho_{r}(t)dt+
\rho_{r}(\alpha)\right)
$$
\noindent where $w(u)$ is the Buchstab function, $\rho_{r}(u)$ is a function which for $0<u\leq1$ equals
$$
\rho_{r}(u)=\frac{u^{r-1}}{\Gamma(r)}
$$
\noindent and for $u>1$ satisfies the differential equation
$$
u\rho_{r}^{\prime}(u)+(1-r)\rho_{r}(u)+r\rho_{r}(u-1)=0,
$$
\noindent and}
$$
C(r)=\prod_{p}\left(1-\frac{1}{p}\right)^r\left(1+\frac{r}{p-1}\right).
$$

\bigskip

\indent The function $\rho_{z}$ arises in the asymptotic analysis of
$$
\Psi_{z}(x,y)=\sum_{\substack{n\leq x\\
		P^{+}(n)\leq y}}z^{\omega(n)}
$$
\noindent where
$$
\Psi_{z}(x,y)\sim\rho_{z}(\alpha)\frac{x}{\log^{1-z}(y)}.
$$
\noindent In particular, $\rho_{z}(\alpha)=e^{-\alpha\log\alpha+O(\alpha)}$.  This rapid decay, coupled with the fact that
$$
\int_{0}^{\infty}\rho_{z}(\alpha)d\alpha=e^{z\gamma},
$$
\noindent ensures that $m_{r}(\alpha)\rightarrow\ell(r)$.

\indent The Selberg method can now be applied to the results of Theorem 1.  The consequence is the following theorem, which provides an estimate of $N_{k}(x,y)$ for \textit{small} values of $y$.

\bigskip

\noindent\textbf{Theorem 10:} \textit{There exists a constant $C>0$ such that if $\alpha>C\log\log(x)$, $k\geq1$, and $r>0$, then
$$
N_{k}(x,y)=\ell\left(\frac{k}{\log\log(y)}\right)\frac{x}{\log(y)}\frac{(\log\log(y))^{k}}{k!}\left(1+O\left(\frac{k}{(\log\log(y))^{2}}\right)\right),
$$
\noindent uniformly for $k\leq r\log\log(y)$.}

\bigskip

\noindent\textbf{Remarks:} Our proof is based upon the method of Selberg
in [10] (generalized in Chapter II.6.1 of [11]).  However, there is an
important difference between his approach and the following proof of
Theorem 10.  In Selberg's estimate for the sum $S_{z}(x)$, there is a
$\frac{1}{\Gamma(z)}$ factor which will, by virtue of the functional
equation of the $\Gamma$-function, absorb the the factor $\frac{1}{z}$  in
$\frac{1}{z^{k+1}}$ in equation (52). In view of this, Selberg's choice for the
radius of the circle in (52), namely
$r=\frac{k-1}{\log\log(x)}$, is optimal due to the vanishing of a first
order error term similar to what is given in (59) below. In contrast,
our estimate of
$S_{z}(x,y)$ does not contain a factor involving the $\Gamma$-function,
and since the function $\ell(z)$ in Theorem 7 has the property that
$\ell(0)\neq0$, it does not absorb any of the factors of $z^{-k-1}$ at $z=0$.
We will see in the course of the proof of Theorem 10 that the optimal choice
for the radius in (52) will be $r=\frac{k}{\log\log(y)}$ (to ensure the
vanishing of the first order error
term in the estimation of $S_{z}(x,y)$).  The absence of the $\Gamma$-function
is the reason why $N_k(x,y)$ is to be compared with $N_{k+1}(x)$ and as will be
seen below. 

\bigskip

\begin{proof} From Theorem 1, with $|z|\leq r$, we have, for 
$$
y\leq x^{1/(r+D+2)\log\log(x)},
$$
\noindent ($D$ the constant in Theorem 1) equivalently, $\alpha\geq(r+D+2)\log\log(x)$
$$
S_{z}(x,y)=x\prod_{p<y}\left(1+\frac{z-1}{p}\right)+O(xe^{-\alpha}\log^{D}(x))+O\left(\frac{x}{\log^{r+2}(x)}\right).
$$
\noindent The above equation is then
$$
S_{z}(x,y)=\ell(z)\frac{x}{\log^{1-z}(y)}+O\left(\frac{x}{\log^{2-z}(y)}\right)
$$
\noindent and so
$$
N_{k}(x,y)=\frac{1}{2\pi i}\int_{|z|=r}S_{z}(x,y)\frac{dz}{z^{k+1}}
$$
$$
=\frac{x}{2\pi i\log(y)}\int_{|z|=r}\frac{\ell(z)\log^{z}(y)}{z^{k+1}}dz+O\left(\frac{x}{\log^{2}(y)}\int_{|z|=r}\frac{\log^{\Re(z)}(y)}{|z|^{k+1}}|dz|\right)
$$
\begin{equation}
=\frac{x}{\log(y)}\frac{1}{2\pi i}\int_{|z|=r}\frac{\ell(z)\log^{z}(y)}{z^{k+1}}dz+O\left(\frac{x}{\log^{2}(y)}\frac{e^{r\log\log(y)}}{r^{k+1}}r\right)
\end{equation}
\noindent as $\Re(z)\leq |z|= r$.  Upon choosing $r=\frac{k}{\log\log(y)}$ equation (53) yields
$$
=\frac{x}{\log(y)}\frac{1}{2\pi i}\int_{|z|=r}\frac{\ell(z)\log^{z}(y)}{z^{k+1}}dz+O\left(\frac{x}{\log^{2}(y)}\frac{e^{k}}{k^{k}}(\log\log(y))^{k}\right)
$$
\begin{equation}
=\frac{x}{\log(y)}\frac{1}{2\pi i}\int_{|z|=r}\frac{\ell(z)\log^{z}(y)}{z^{k+1}}dz+O\left(\frac{x}{\log^{2}(y)}\frac{\sqrt{k}(\log\log(y))^{k}}{k!}\right)
\end{equation}
\noindent from a weak form of Stirling's formula $k!<<k^{k+1/2}e^{-k}$.
	
\indent Recall Theorem 7 where it was noted that $\ell(z)$ is analytic, and so we may write 
\begin{equation}
\ell(z)=\ell(r)+\ell^{\prime}(r)(z-r)+O((z-r)^{2})
\end{equation}
\noindent so that using (55) we may represent the integral on the right of (54) in the following form
$$
\frac{x}{\log(y)}\frac{1}{2\pi i}\int_{|z|=r}\frac{\ell(z)\log^{z}(y)}{z^{k+1}}dz
$$
$$
=\frac{x}{\log(y)}\frac{\ell(r)}{2\pi i}\int_{|z|=r}\frac{e^{z\log\log(y)}}{z^{k+1}}dz+\frac{x}{\log(y)}\frac{\ell^{\prime}(r)}{2\pi i}\int_{|z|=r}(z-r)\frac{e^{z\log\log(y)}}{z^{k+1}}dz
$$
$$
+O\left(\frac{x}{\log(y)}\int_{|z|=r}(z-r)^{2}\frac{\log^{\Re(z)}(y)}{z^{k+1}}dz\right)	$$
$$
=I_{1}+I_{2}+I_{3}.
$$
\noindent The integral $I_{1}$ may be easily evaluated by direct appeal to the Cauchy integral formula:
\begin{equation}	I_{1}=\frac{x}{\log(y)}\frac{\ell(r)}{2\pi i}\int_{|z|=r}\frac{e^{z\log\log(y)}}{z^{k+1}}dz=\ell\left(\frac{k}{\log\log(y)}\right)\frac{x}{\log(y)}\frac{(\log\log(y))^{k}}{k!}
\end{equation}
\noindent upon identifying the integral with the coefficient of $z^{k}$ in the series expansion of $e^{z\log\log(y)}$.
\indent With regard to the integral $I_{2}$, note that
$$
\frac{\ell^{\prime}(r)}{2\pi i}\int_{|z|=r}(z-r)\frac{e^{z\log\log(y)}}{z^{k+1}}dz
$$
$$
=\frac{\ell^{\prime}(r)}{2\pi i}\int_{|z|=r}\frac{e^{z\log\log(y)}}{z^{k}}dz-\frac{\ell^{\prime}(r)}{2\pi i}r\int_{|z|=r}\frac{e^{z\log\log(y)}}{z^{k+1}}dz
$$
\begin{equation}
=\ell^{\prime}(r)\frac{(\log\log(y))^{k-1}}{(k-1)!}-\ell^{\prime}(r)r\frac{(\log\log(y))^{k}}{k!}.
\end{equation}
\noindent With the choice of $r=\frac{k}{\log\log(y)}$, equation (57) becomes
\begin{equation}
I_{2}=\ell^{\prime}(r)\frac{(\log\log(y))^{k-1}}{(k-1)!}-\ell^{\prime}(r)\frac{(\log\log(y))^{k-1}}{(k-1)!}=0.
\end{equation}
\noindent Equation (58) demonstrates why choosing $r=\frac{k}{\log\log(y)}$ is, essentially, the best possible choice.  Any other value for $r$ in (58) would contribute an error term which would not be zero.
	
\indent It remains to estimate the third integral $I_{3}$.  Let us set $z=re^{2\pi i\theta}$ for $-\frac{1}{2}\leq\theta\leq\frac{1}{2}$, so that 
$$
z-r=2ie^{\pi i\theta}r\frac{(e^{\pi i\theta}-e^{-\pi i\theta})}{2i}=2ie^{\pi i}r\sin\pi\theta<<r\theta.
$$
\noindent Thus
$$
(z-r)^{2}<<r^{2}\theta^{2}.
$$
\noindent Next note that
$$
|\log^{z}(y)|=|e^{z\log\log(y)}|=e^{\Re(z)\log\log(y)}=e^{r\cos2\pi\theta\log\log(y)}	$$
\begin{equation}
\leq e^{r(1-\lambda\theta^{2})\log\log(y)}=e^{k(1-\lambda\theta^{2})}
\end{equation}
\noindent by again using the fact that $r=k/\log\log(y)$, and where $\lambda>0$ is some small fixed constant.  Applying the bound in (59) to the integral $I_{3}$ we see
$$
I_{3}<<\int_{-1/2}^{1/2}\frac{xr^{2}\theta^{2}}{\log(y)}\frac{e^{k}e^{-k\lambda\theta^{2}}r|e^{2\pi i\theta}|}{r^{k+1}}d\theta
$$
$$
<<\frac{x}{\log(y)}\frac{e^{k}}{r^{k-2}}\int_{-1/2}^{1/2}\theta^{2}e^{-k\lambda\theta^{2}}d\theta<<\frac{x}{\log(y)}\frac{e^{k}}{r^{k-2}}\int_{0}^{\infty}\theta^{2}e^{k\lambda\theta^{2}}d\theta
$$
\noindent which after the substitution $\sqrt{k}\theta=u$ yields
$$
\frac{x}{\log(y)}\frac{e^{k}}{r^{k-2}}\int_{0}^{\infty}\frac{u^{2}e^{-\lambda u^{2}}}{k\sqrt{k}}du
$$
$$
<<\frac{x}{\log(y)}\frac{(\log\log(y))^{k-2}}{k^{k-2}}\frac{e^{k}}{k^{3/2}}<<\frac{x}{\log(y)}\frac{(\log\log(y))^{k-2}}{e^{-k}k^{k-1/2}}
$$
\begin{equation}
<<\frac{x}{\log(y)}\frac{(\log\log(y))^{k}}{k!}\frac{k}{(\log\log(y))^{2}},
\end{equation}
\noindent by using Stirling's formula once more.  Combining (60) with the estimates for $I_{1}$ and $I_{2}$ gives
$$
N_{k}(x,y)=I_{1}+I_{2}+I_{3}=\ell(r)\frac{x}{\log(y)}\frac{(\log\log(y))^{k}}{k!}\left(1+O\left(\frac{k}{(\log\log(y))^{2}}\right)\right)
$$
$$
=\ell\left(\frac{k}{\log\log(y)}\right)\frac{x}{\log(y)}\frac{(\log\log(y))^{k}}{k!}\left(1+O\left(\frac{k}{(\log\log(y))^{2}}\right)\right).
$$
\end{proof}

\indent It is to be noted that for almost all integers we have $\omega_{y}(n)\sim\log\log(y)$, and for this situation
$$
\frac{x}{\log(y)}\frac{(\log\log(y))^{k}}{k!}\sim\frac{x}{\log(x)}\frac{(\log\log(x))^{k-1}}{(k-1)!}.
$$
\noindent On the other hand, if $k\sim\lambda\log\log(y)$ with $\lambda\neq1$ then
$$
\frac{x}{\log(y)}\frac{(\log\log(y))^{k-1}}{(k-1)!}\sim\lambda\frac{x}{\log(y)}\frac{(\log\log(y))^{k}}{k!},
$$
\noindent that is, the terms differ by a factor of $\lambda$.

\indent We will now study $N_{k}(x,y)$ for \textit{large} $y$.  Tenenbaum [14]
has communicated to us that by suitably adapting the powerful techniques of
[13] to $S_{z}(x,y)$, one may derive an effective estimate for $N_{k}(x,y)$
for certain ranges of $k$.  These estimates follow from a suitable upper
bound for $S_{z}(x,y)$ when $|z|=r$ but is not close to $r$, and asymptotic
estimates for $S_{z}(x,y)$ when $z$ is close to $r$.  The form in which we
will use his result is

\bigskip

\noindent\textbf{Theorem 11:} \textit{Let $\kappa>0$ be a small parameter,
$r=\frac{k}{\log\log(y)+c_{1}}$, where
$$
c_{1}:=\lim_{y\rightarrow\infty}\sum_{p<y}\frac{1}{p}-\log\log(y).
$$
\noindent Then for $\kappa\leq r\leq1/\kappa$,
$$
N_{k}(x,y)=S_{r}(x,y)\frac{(\log\log(y)+c_{1})^{k}}{k!e^{k}}\left(1+O\left(\frac{1}{\sqrt{\log\log(y)}}\right)\right)
$$
\noindent uniformly for $1\leq\alpha<(\log\log(x))^{2}$.}

\bigskip

\noindent\textbf{Remark:} If one chooses $r=\frac{k}{\log\log(y)}$ in Theorem 11 then the result becomes
\begin{equation}
N_{k}(x,y)=S_{r}(x,y)\frac{(\log\log(y))^{k}}{k!e^{k}}\left(1+O\left(\frac{1}{\sqrt{\log\log(y)}}\right)\right)
\end{equation}
\noindent which is more convenient.

\indent Theorem 11 interesting for several reasons. Firstly, it
gives a relationship between the coefficients $N_{k}(x,y)$ of
$z^k$ in the sum $S_{z}(x,y)$ with the sum itself.  Secondly, the size of
sum $S_{z}(x,y)$ is utilized only when $z=r$ is real-valued and
positive.  Thankfully, the results of Section 4 apply in this situation
and so we get the following theorem.

\bigskip

\noindent\textbf{Theorem 12:} \textit{Let $\kappa>0$, $r=\frac{k}{\log\log(y)}$, and $\kappa<r<1/\kappa$, then
$$
N_{k}(x,y)= m_{r}(\alpha)\frac{x(\log\log(y))^{k}}{k!\log(y)}\left(1+O\left(\frac{1}{\sqrt{\log\log(y)}}\right)\right)
$$
\noindent uniformly for $1\leq\alpha\leq(\log\log(x))^{2}$.}

\bigskip

\begin{proof} \noindent Applying the estimate in Theorem 6 this with $\alpha<(\log\log(x))^{2}$ to equation (61) yields
\begin{equation}
N_{k}(x,y)= m_{r}(\alpha)\frac{x(\log\log(y))^{k}}{k!e^{k}\log^{1-r}(y)}\left(1+O\left(\frac{1}{\sqrt{\log\log(y)}}\right)\right).
\end{equation}
\noindent As $k=r\log\log(y)$ (62) can be rewritten
	
$$
N_{k}(x,y)=m_{r}(\alpha)\frac{x(\log\log(y))^{k}}{k!\log(y)}\left(1+O\left(\frac{1}{\sqrt{\log\log(y)}}\right)\right).
$$
	
\end{proof}

\noindent\textbf{Remark:} We saw in section 4 that when estimating $S_{z}(x,y)$ for \textit{large} $y$ we were forced to restrict $\Re(z)>0$, so that a direct application of the Selberg contour integral method would not be possible without estimates for the case $\Re(s)\leq0$.  However, with suitable bounds for $\Re(z)\leq0$ the method can still yield the correct asymptotic equality.  We note that the methods employed in Section 4 will provide us with the bound of

$$
|S_{z}(x,y)|<<\frac{x(\log\log(x))^{K+1}}{\log(y)},
$$  

\noindent which is uniform provided $\alpha<<\log\log(x)$, $\Re(z)\leq0$, and $|z|\leq R$ (and where $K$ is the constant in Theorem 6).  Note that this bound holds in the range $\alpha<<\log\log(x)$ which overlaps with the range of $\alpha$ in section 3.  Therefore, we could have used this bound (which is milder than the results of Theorem 11).  We made use of Theorem 11 to conclude the main result of Theorem 12 because of the sharpness of the error term in this theorem.

\indent We note that when $k$ is fixed, then $k/\log\log(y)$ is close to $0$.  But then $m_{r}(\alpha)\rightarrow m_{0}(\alpha)=w(\alpha)$, which means that Theorem 12 in this case corresponds asymptotically to Corollary 2. Also, if $y=x$
then Theorem 12 implies
$$
N_{k}(x,y)\sim\frac{S_{r}(x)}{e^{k}}\frac{(\log\log(y))^{k}}{k!}=\frac{g(1,r)}{\Gamma(r)}\frac{x}{\log(y)}\frac{(\log\log(y))^{k}}{k!}
$$
$$
=\frac{g(1,r)r}{\Gamma(1+r)}\frac{x}{\log(x)}\frac{(\log\log(x))^{k}}{k!}=\frac{g(1,r)}{\Gamma(1+r)}\frac{x}{\log(x)}\frac{(\log\log(x))^{k-1}}{(k-1)!},
$$
\noindent which corresponds to Selberg's estimate (3) for $N_{k}(x)$.

\indent It is a matter of taste how one presents the estimate for
$N_{k}(x,y)$.  In Theorem 12 we wrote this in terms of ratios of elementary
functions, and in the context of the results of the previous section, this
is the more natural way to present this estimate.  Tenenbaum has communicated
Theorem 11 to us as a special case of his general result in [13] on the local
distribution of $\omega(n;E)$.  However, as $E$ is a general set of primes,
it is difficult to expect asymptotic estimates in terms of elementary
continuous functions.  Thus, Tenenbaum's result (Cor 2.4 of [13]) gives an
estimate for the local distribution of $\omega(n;E)$ in terms of
$S_{r}(x,E)$.  If $E$ has regular behavior, such as when $E=\{p|p<y\}$,
then $S_{r}(x,E)=S_{r}(x,y)$ can be estimated in terms of continuous functions.

\bigskip

\noindent \textit{ACKNOWLEDGMENTS:} The results presented here are contained
in the PhD thesis of Todd Molnar submitted to the University of Florida in
2017 under
the guidance of Krishnaswami Alladi. We are thankful to Professor Gerald
Tenenbaum who generously shared his results in a letter. 

\section{Bibliography}
\noindent [1] K. Alladi, \textit{Asymptotic Estimates of Sums involving the Moebius Function}, J. Number Theory \textbf{14} (1982), 86-98.

\medspace

\noindent [2] K. Alladi, \textit{Asymptotic Estimates of Sums Involving the Moebius Function. II}, Trans. Amer. Math. Soc. \textbf{272}, (1982), 87-105.

\medspace

\noindent [3] K. Alladi, \textit{The distribution of $\nu(n)$ in the sieve of Eratosthenes}, Quart. J. Math. Oxford (2), \textbf{33} (1982), 129-148.

\medspace

\noindent [4] K. Alladi, \textit{Multiplicative functions and Brun's Sieve}, Acta Arith. \textbf{51}, 201-219 (1988).

\medspace

\noindent [5] N. G. de Bruijn, \textit{On the number of uncanceled elements in the sieve of Eratosthenes}, Indag. Math. \textbf{12} (1950), 247-256.

\medspace

\noindent [6] P.D.T.A. Elliott, \textit{Probabilistic Number Theory, vol. 1 and 2}, Grundlehren vol. 239-240 Springer-Verlag, Berlin and New York, 1979.

\medspace

\noindent [7] G. Hal\'{a}sz, \textit{Remarks to my paper 'On the distribution of additive and the mean values of multiplicative arithmetic functions'}, Acta Math. Acad. Scient. Hung. \textbf{23}, 425-432.

\medspace

\noindent [8] R.R. Hall and G. Tenenbaum, \textit{Divisors}, Cambridge tracts in mathematics, no. 90, Cambridge University Press, 1988.

\medspace

\noindent [9] L. G. Sathe, \textit{On a problem of Hardy on the distribution of integers with a given number of prime factors I-IV}, J. Indian Math. Soc. (N.S.) \textbf{17}, 63-82 and 83-141, \textbf{18} 27-42 and 43-81.

\medspace

\noindent [10] A. Selberg, \textit{Note on a paper of L.G. Sathe}, J. Indian Math. Soc. \textbf{18} (1954), 83-87.

\medspace

\noindent [11] G. Tenenbaum, \textit{Introduction to analytic and probabilistic number theory}, 3rd ed., Graduate Studies in Mathematics {\textbf{163}},
Amer. Math. Soc. 2015.

\medspace

\noindent [12] G. Tenenbaum, \textit{Generalized Mertens sums}, in Analytic
Number Theory, Modular Forms, and $q$-Hypergeometric Series - in honor of
Krishna Alladi's 60th Birthday (G. E. Andrews and F. Garvan, Eds.), Springer
Proc. in Math. and Stat. {\textbf{221}} (2018), 733-736.

\medspace

\noindent [13] G. Tenenbaum, \textit{Valeurs moyennes effectives de fonctions
multiplicatives complexes}, Ramanujan J. \textbf{44} (2017), 641-701.

\medspace

\noindent [14] G. Tenenbaum, \textit{Private Communication} (2016).

\bigskip
\bigskip

Department of Mathematics

University of Florida

Gainesville, FL 32611, USA

{\textit{email:}} alladik(@ufl.edu

{\textit{email:}} twmolnar@gmail.com

\end{document}